\documentclass{article}
\usepackage{diagrams}

\usepackage{mathptmx}

\usepackage[textwidth=12.8cm, textheight=18.5cm]{geometry}

\usepackage{amssymb}
\usepackage{amsmath}

\newtheorem{prop}{Proposition}
\newtheorem{thm}[prop]{Theorem}
\newtheorem{corol}[prop]{Corollary}

\title{Commutative monads as a theory of distributions}
\author{Anders Kock\\
University of Aarhus}
\date{}

\DeclareMathOperator{\tot}{tot}

\newcommand{\A}{\mathcal{A}}

\newcommand{\E}{\mathcal{E}}
\newcommand{\p}{\pitchfork}

\newcommand {\imes}{\times}
\newcommand{\ev}{ev}
\begin{document}
\maketitle
\section*{Introduction}
The word ``distribution'' has a common, non technical, meaning, 
roughly synonymous with ``dispersion'': a given  quantity, say  
a quantity of 
tomato, may be distributed or dispersed over a given space, say
a pizza. A more mathematical use of the word  
was made precise in functional analysis by L.\ Schwartz (and 
predecessors):  it is a 
(continuous linear) functional on a space of functions.
Schwartz argues that ``the mathematical distributions constitute a 
correct mathematical definition of the distributions one meets in 
physics'', \cite{Schwartz} p.\ 84.

Our aim here is to present an alternative mathematical theory of 
distributions (of compact support), applicable also for the ``tomato'' 
example, but which does not depend on the ``double dualization'' 
construction of functional analysis; and also, to provide a canonical 
comparison to the distributions of functional analysis.

Distributions of compact  support form an important example of an {\em 
extensive quantity}, in the sense made precise by Lawvere, cf.\ e.g.\ 
\cite{SQ}, and whose salient feature is   the covariant functorial 
dependence 
on the space over which it is distributed. 
Thus, there is a covariant functor $T$, such that 
extensive quantities of a given type on a given space $X$ form a new 
``linear'' space 
$T(X)$. In the theory we present here, $T$ is an endofunctor on a 
cartesian closed category $\E$, in fact a monad; a linear 
structure grows, for monads with a certain property, out of the monad structure.  

An example of such $T$ is, for suitable $\E$, the double dualization 
construction $S$ that gives the space of Schwartz distributions of 
compact support, -- except that this $S$ is not commutative. The unit 
$\eta$ 
of the monad $S$ associates to an $x\in X$ the Dirac distribution 
$\delta_{x}$ at $x$.

What makes our theory simple is the universal property which the unit 
map $\eta : X\to T(X)$ is known to have by 
general monad theory, -- in conjunction with the assumed commutativity 
of $T$ (in the sense of \cite{MSMCC}). This is what makes  for 
instance the notions of
linearity/bilinearity work well.

\medskip
\noindent{\bf Generalities:} We  mainly  compose maps from right to 
left (this is the default, and it is denoted $g\circ f$); but occasionally, in 
particular in connection with displayed diagrams, we compose from left to right 
(denoted $f.g$); we remind the reader when this is the case.

In Sections \ref{ad1x}-\ref{ad12x}, the notions and theory (except 
the remarks on ``totals'') work for 
$\E$ a general symmetric monoidal closed category, not necessarily 
cartesian closed. Here, the reasoning is mainly diagrammatic. In the 
remainder, we take the liberty to reason with ``elements''. In 
particular, the ``elements'' of $T(X)$ are simply called 
``distributions on $X$''

\medskip

The present text subsumes and simplifies the preliminary arXiv texts, 
\cite{MEQ}, \cite{CEQ}, and has been presented in part in Krakow 
at the conference ``Differential Geometry and Mathematical Physics, June-July 2011,
in honour of Wlodzimierz Tulczyjew'', 
 and at the Nancy Symposium ``Sets within Geometry'' July 2011. 
I want to thank the organizers of these meetings 
for inviting me. I also want to thank 
 Bill Lawvere for fruitful discussions in Nancy, 
and  e-mail correspondence in the spring of 2011.

\section{Monads, algebras, and linearity}\label{ad1x}
Recall the notion of monad $T=(T,\eta,\mu )$ on a category $\E$, cf.\ 
e.g.\ \cite{Mac} 6.1 and 6.2. Recall also the notion of $T$-algebra 
$A=(A,\alpha )$ for such monad; here, $\alpha :T(A)\to A$ is the 
{\em structure} map for the given algebra. There is a notion of {\em 
morphism} of algebras $(A,\alpha )\to (B, \beta )$, cf.\ loc.cit., so 
that we have the category $\E ^{T}$ of algebras for 
the monad $T$. For each $X\in \E$, we have the $T$-algebra $T(X)$ 
with structure map $\mu_{X}:T^{2}(X)\to T(X)$. The map $\eta 
_{X}:X\to T(X)$ has a universal property, making $(T(X),\mu _{X})$ 
into a {\em free} algebra on $X$: to every $T$-algebra $(A,\alpha )$ 
and every map $f:X \to A$ in $\E$, there is a unique morphism of 
algebras $\overline{f}: T(X)\to A$ with
$\overline{f} \circ \eta _{X}= f$ (namely $\overline{f}=\alpha \circ 
T(f)$).

Recall that any algebra structure $\alpha$ on an object $A$ is by 
itself an algebra morphism $T(A) \to A$, and in particular $\mu _{X}$ 
is an algebra morphism. Recall also that any morphism $T(f): T(X)\to 
T(Y)$ (for $f:X\to Y$ an arbitrary morphism in $\E$) is an algebra 
morphism.

For the monads that we are to consider in the present article, namely  
{\em commutative} monads (cf.\ Section \ref{ad10x} below), there are good reasons for 
an alternative  terminology: namely, $T$-algebras deserve the name 
{\em $T$-linear spaces}, and homomorphisms deserve the name 
{\em $T$-linear maps} (and, if $T$ is understood from the context, 
the `$T$' may 
even be omitted). This allows us to talk about {\em partial $T$-linear 
maps}, as well as {\em $T$-bilinear maps}, as we shall explain. 

An example of a commutative monad, with $\E$ the category of sets, is 
the functor $T$ which to a set $X$ associates (the underlying set of) 
the free real vector space on $X$. In this case, the algebras for $T$ 
are the  vector spaces over ${\mathbb R}$, and the morphisms are the 
${\mathbb R}$-linear maps.

For a general monad, the ``linear'' terminology is usually not 
justified, but we shall nevertheless use this terminology right from 
the beginning, in so far as maps are concerned. Thus, the map $\overline{f}:T(X) \to A$ considered 
above will be called {\em the $T$-linear extension of $f: X\to A$} 
(and $\eta_{X}:X\to T(X)$ is then to be understood from the context). 

\section{Enrichment and strength}\label{ad2x}
We henceforth consider the case where $\E$ is a cartesian closed 
category (cf.\ e.g.\ \cite{Mac} 4.6); so, for any pair of objects 
$X$, $Y$ in $\E$, we have the exponential object $Y^{X}$, and an 
``evaluation'' map $\ev : Y^{X}\times X \to Y$ (or more pedantically: 
$\ev _{X,Y}$). Since we shall 
iterate the ``exponent formation'', it becomes expedient to have an 
on-line notation. Several such are in use in the literature, like $[X,Y]$, $X\multimap 
Y$, or $X\p Y$; we shall use the latter (tex code for $\p$ is \verb+ \pitchfork+; 
read ``$X\p Y$'' as ``$X$ hom $Y$''), so that the evaluation map is a map
$$\ev : (X\p Y) \times X \to Y.$$
Since $\E$ is  a  monoidal category (with cartesian product as 
monoidal structure), one has the notion of 
when a category is {\em enriched} in $\E$, cf.\ e.g.\ \cite{Borceux} 
II.6.2. And since this monoidal category $\E$ is in fact monoidal 
closed (with $\p$ as the closed structure), $\E$ is enriched in 
itself. For $\E$-enriched categories, one has the notion of enriched 
functor, as well as enriched natural transformation between such. 
(cf.\ e.g.\ loc.cit.\ Def.\ 6.2.3 and 6.2.4). So in particular, it 
makes sense to ask for an $\E$-enrichment of the functor $T:\E \to 
\E$ (and also to ask for $\E$-enrichment of $\eta$ and $\mu$, which 
we shall, however, not consider until in the next Section).

\begin{sloppypar}Specializing the definition of enrichment to the case of an 
endofunctor $T:\E \to \E$ gives that such enrichment consists in maps 
in $\E$
$$\begin{diagram}X\p Y &\rTo ^{st_{X,Y}}& T(X)\p T(Y)\end{diagram},$$
(for any pair of objects $X,Y$ in $\E$)
satisfying a composition axiom and a unit axiom.
The word ``enriched'' is sometimes replaced by the word ``strong'', 
and ``enrichment'' by ``strength'', whence the notation ``$st$''. We 
shall, however, need to consider two other equivalent mani\-festations 
of such ``strength'' 
structure on $T$, introduced in \cite{MSMCC} (and proved equivalent 
to strength in \cite{SFMM}, Theorem 1.3), and in  \cite{CCGBCM}, 
called {\em tensorial} and {\em cotensorial} strength, respectively.
To give a tensorial strength to the endofunctor $T$ is to give, for 
any pair of objects $X,Y$ in $\E$ a map
$$\begin{diagram}X\times T(Y) &\rTo^{t''_{X,Y}}& T(X\times 
Y)\end{diagram},$$ satisfying a couple of 
equations (cf.\ \cite{SFMM} (1.7) and (1.8)). The tensorial strength $t''$ has a ``twin sister'' 
$$\begin{diagram}T(X)\times Y &\rTo^{t'_{X,Y}}& T(X\times 
Y)\end{diagram},$$ essentially obtained by 
``conjugation $t''$ by the twist map $X\times Y \to Y\times X)$''. Similarly, to give a cotensorial strength to the 
endofunctor $T$
 is to give, for 
any pair of objects $X,Y$ in $\E$ a map
$$\begin{diagram}T(X\p Y) &\rTo^{\lambda_{X,Y}}&X\p T(Y)\end{diagram},$$
(cf.\ \cite{CCGBCM}) satisfying a couple of 
equations.  
The three manifestations of strength, and also the equations they 
have to satisfy, are deduced from one another by simple ``exponential 
adjointness''-transpositions, cf.\ loc.cit. As an example, let us 
derive the tensorial strength  $t''$ from the classical ``enrichment'' strength 
$st$:
We have the unit $u$ of the adjunction $ (-\times Y)\dashv (Y\p -)$; it 
is the first map in the following composite; the second is the 
enrichment strength
$$\begin{diagram}X&\rTo^{u}&Y\p (X\times Y)&\rTo^{st_{Y,X\p Y}}& 
T(Y)\p T(X\times Y). 
\end{diagram}$$
Now apply exponential adjointness to transform this composite map to the desired $t''_{X,Y}: 
X\imes T(Y) \to T(X\times Y)$.\end{sloppypar}

The strength in its tensorial form $t''$ (and its twin sister $t'$) 
will be the main manifestation used in the present paper. The 
tensorial form has the advantage not of mentioning $\p$ at all, so 
makes sense even for endofunctors on monoidal categories that are not 
closed; this has been exploited e.g.\ in \cite{Moggi} and \cite{CJ}. 

\section{Strong monads; partial linearity}\label{ad3x}
The composite of two strong endofunctors on the cartesian closed $\E$ 
carries a strength derived from the two given strengths. We give it 
explicitly for the composite $T\circ T$ (with strength of $T$ given in 
tensorial form $t''$), namely it is the ``composite'' strength is the 
top line in the diagram (\ref{multipurpose}) below. 

There is also an evident notion of strong natural transformation between 
two strong functors;  in terms of tensorial strength, and for the 
special case of the natural transformation $\mu :T\circ T \Rightarrow 
T$, this means that (not only  the naturality squares but also)
 the top square in (\ref{multipurpose}) commutes, for all $X,Y$.

The identity functor $I: \E \to \E$ has identity maps $id_{X}$ as its 
tensorial strength. The bottom square in 
(\ref{multipurpose}) expresses that the natural transformation $\eta 
:I \Rightarrow T$ is strongly natural.

The notion of composite strength, and of strong natural transformation, 
are equivalent to the classical notion of strength/enrichment (in 
terms of $st$), in so 
far as endofunctors on $\E$ goes, see \cite{MSMCC} (and more 
generally, for functors between categories enriched over $\E$, 
provided they are tensored over $\E$, cf.\ \cite{SFMM}).

Here is the diagram which expresses the strength of $\mu :T^{2} 
\Rightarrow T$ (upper part), and of $\eta :I \Rightarrow T$ (lower part):

\begin{equation}\label{multipurpose}\begin{diagram}X\times 
T^{2}Y&\rTo^{t''_{X,TY}}&T(X\imes 
TY)&\rTo^{T(t''_{X,Y})}&T^{2}(X\times Y)\\
\dTo^{X\times \mu_{Y}}&&&&\dTo_{\mu_{X\times Y}}\\
X\times TY&&\rTo_{t''_{X,Y}}&&T(X\times Y)\\
\uTo^{X\times \eta _{Y}}&&&&\uTo_{\eta_{X\times Y}}\\
X\times Y&&\rTo_{id}&&X\times Y
\end{diagram}
\end{equation}

One can equivalently formulate strength of composite functors, and 
strength of natural transformations, in terms of  cotensorial 
strength, cf.\ \cite{CCGBCM}; this will not explicitly play 
a role here. But one consequence of the cotensorial strength will be 
important: if $(C, \gamma)$ is an algebra for the monad $T$, and $X$ 
is any object in $\E$, then $X\p C$ carries a canonical structure of 
algebra for $T$, with structure map the composite
$$\begin{diagram}T(X\p C) &\rTo^{\lambda _{X,C}}&X\p T(C)&\rTo^{X\p 
\gamma }&X\p C\end{diagram}.$$
The $T$-algebra $X\p C$ thus constructed actually witnesses that the 
category $\E ^{T}$ of $T$-algebras is {\em cotensored} over $\E$, 
cf.\ e.g.\ \cite{Borceux} II.6.5.


\medskip
The tensorial strength makes possible a description (cf.\ \cite{BCCM}) of the crucial notion of 
{\em partial $T$-linearity} (recall the use of the phrase {\em $T$-linear} 
as synonymous with $T$-algebra morphism). Let
$(B,\beta )$ and $(C, \gamma )$ be $T$-algebras, and let $X \in E$ be an 
arbitrary object. Then a map $f:X\times B \to C$ is called {\em $T$-linear 
in the second variable} (or {\em 2-linear} , if $T$ can be understood 
from the context) if the following pentagon commutes:
\begin{equation}\label{2lin}\begin{diagram}
X\times T(B)&\rTo^{t''_{X,B}}&T(X\times B)&\rTo^{T(f)}&T(C)\\
\dTo^{X\times \beta}&&&&\dTo_{\gamma}\\
X\times B &&\rTo_{f}&&C.
\end{diagram}\end{equation}
In completely analogous way, one describes the notion of 
{1-linearity} of a map $A\times Y \to C$, where $A$ and $C$ are 
equipped with algebra structures: just apply $t'_{A,Y}$. Finally, a 
map $A\times B \to C$ is called {\em bilinear} if it is both 1-linear 
and 2-linear.

The notion of 2-linearity of $f: X\times B \to C$ can equivalently, 
and perhaps more intuitively, be described:
the exponential transpose of $f$, i.e.\ the map $\hat{f}:B\to X\p C$, 
is $T$-linear (= a morphism of $T$-algebras), with $T$-structure on 
$X\p C$ as given above in terms of the cotensorial strength $\lambda$.

We shall prove
\begin{prop}\label{onexx} The map $t''_{X,Y}:X\times T(Y)\to 
T(X\times Y)$ is 2-linear; and it is in fact initial among 2-linear 
maps from $X\times T(Y)$ into $T$-algebras. (Similarly, $t'_{X,Y}$ is 
1-linear, and is initial among 1-linear maps from $T(X)\times Y$ into 
$T$-algebras.)
\end{prop}
{\bf Proof.} We first argue that $t''_{X,Y}$ is 2-linear. This means 
by definition that a certain pentagon commutes; for the case of 
$t''_{X,Y}$, this is precisely the (upper part of) the diagram 
(\ref{multipurpose}) 
expressing that $\mu$ is strongly natural. 

To prove that $t''_{X,Y}$ is initial among 2-linear maps: Given a 
2-linear $f:X\times T(Y) \to B$,
where $B=(B,\beta)$ is a $T$-algebra. We must prove 
unique existence of a linear $\overline{f}$ making the right hand 
triangle in the following diagram commutative:
$$\begin{diagram}[nohug]X\imes Y&\rTo^{\eta_{X\times Y}}&T(X\times Y)\\
\dTo^{X\times \eta_{Y}}&\ruTo_{t''_{X,Y}}&\dTo_{\overline{f}}\\
X\imes TY&\rTo_{f}&B.
\end{diagram}$$
Since the upper triangle commutes (strong naturality of $\eta$, cf.\ 
the lower part of (\ref{multipurpose})), we see that $\overline{f}$ 
necessarily is the (unique) linear extension over $\eta_{X\times Y}$ 
of $$\begin{diagram}X\times Y &\rTo ^{X\imes \eta_{Y}}&X\imes 
TY&\rTo^{f}&B\end{diagram},$$
thus $\overline{f}$ is, by standard monad theory, the composite of the three last arrows in
\begin{equation}\label{ofx}\begin{diagram}
X\imes TY&\rTo^{t''_{X,Y}}&T(X\imes Y)&\rTo ^{T(X\times 
\eta_{Y})}&T(X\times TY)&\rTo ^{Tf}&TB&\rTo^{\beta}&B
\end{diagram}.\end{equation}
 To prove that this 
$\overline{f}$ does indeed solve the posed problem, we have to prove 
that the total composite in (\ref{ofx}) equals $f$. Composing forwards, 
this follows from the equations
$$\begin{array}{lcll}t''_{X,Y}.T(X\times \eta _{Y}).Tf. \beta &=&
X\times T(\eta_{Y}).t''_{X,TY}.Tf.\beta&\mbox{\quad by naturality of $t''$ w.r.to $\eta_{Y}$}\\
&=&X\imes T(\eta_{Y}).X\times \mu_{Y}.f&\mbox{\quad  by the assumed 
2-linearity of $f$;}\end{array}$$
but this gives  $f$, by a monad law.

\medskip 

It is a direct consequence of the definitions that postcomposing a 
2-linear map $f:X\times B \to C$ with a linear map again yields a 2-linear map. Likewise, 
precomposing $f$ with a map of the form $X\times g$, with $g$ linear, 
again yields a 2-linear map. Similarly for 1-linearity and 
bilinearity.   
\medskip

\section{Partially linear extensions}\label{ad4x}

 The universal property of $\eta_{X}$ (making $T(X)$ the 
``free $T$-algebra on $X$'') was mentioned in Section \ref{ad1x}. 
There are 
similar universal properties of $X\times \eta _{Y}$ and of $\eta_{X} \imes 
Y$:
\begin{prop}\label{one-x}Let $(B,\beta )$ be a $T$-algebra. To any 
$f:X\times Y \to B$, there exists a unique 2-linear 
${\overline{f}}:X\times T(Y) \to B$ making the triangle 
$$\begin{diagram}[nohug]X\times TY&\rTo^{\overline{f}}&B\\
\uTo^{X\times \eta _{Y}}&\ruTo_{f}&&\\
X\times Y&&&
\end{diagram}$$
commute. (There is a similar universal property of $\eta_{X}\times 
Y:X\times Y \to T(X)\times Y$ for 1-linear maps.)
\end{prop}
{\bf Proof.} In essence, this follows by passing to the exponential 
transpose $\hat{f}:Y\to X\p B$ of $f$,  and then using the universal property of 
$\eta_{Y}$. 
More explicitly, there are 
natural bijective correspondences
$$\hom (X\times Y,B)\cong\hom (X,Y\p B) \cong \hom_{T}(T(X), Y\p B)$$
(where the second occurrence of $Y\p B$ is the cotensor  $Y\p B $ in $\E ^{T}$,
 and where the second bijection is induced by 
precomposition by $\eta _{X}$); finally, the 
set $\hom_{T}(T(X), Y\p B)$ is in bijective correspondence with the 
set of 1-linear maps $T(X)\times Y \to B$, by \cite{BCCM} Proposition 1.3 
(i).

A  direct description  of $\overline{f}$ is
\begin{equation}\label{2-lin-ext}\begin{diagram}X\times TY&\rTo ^{t''_{X,Y}}&T(X\times 
Y)&\rTo^{T(f)}&TB&\rTo^{\beta}&B
\end{diagram}.\end{equation}
For, the displayed map is 2-linear, being a composite of the 2-linear 
$t''_{X,Y}$ and the linear maps $T(f)$ and $\beta$; and its precomposition 
with $X\times \eta _{Y}$ is easily seen to be $f$, using that $\eta$ 
is a strong natural transformation.

Similarly, and explicit formula for the 1-linear extension of 
$f:X\times Y \to B$ is
\begin{equation}\label{1-lin-ext}\begin{diagram}T(X)\times Y&\rTo ^{t'_{X,Y}}&T(X\times 
Y)&\rTo^{T(f)}&TB&\rTo^{\beta}&B
\end{diagram}.\end{equation}

\medskip

A consequence of Proposition \ref{one-x} is that $t''_{X,Y}$ may be 
characterized ``a posteriori'' as the unique 2-linear extension over $X\times 
\eta_{Y}$ of $\eta_{X\times Y}$.  But note that $t''$ by itself is 
described independently of $\eta$ and $\mu$. 
Similarly for $t'$.

Also, $\lambda _{X,Y}:T(X\p Y) \to X\p T(Y)$ may be characterized 
as the unique linear extension over $\eta_{X\p Y}$ of $X\p \eta_{Y}$.

\section{The two ``Fubini'' maps}\label{ad6x} The map $t'_{X,Y}:T(X)\times Y\to 
T(X\times Y)$ extends by Proposition \ref{one-x} uniquely over 
$T(X)\times \eta_{Y}$ to a 2-linear map $T(X)\times T(Y)\to T(X\times 
Y)$, which we shall denote $\otimes$ or $\otimes_{X,Y}$. Thus, 
$\otimes$ is characterized as the unique 2-linear map making the 
following triangle commute:
\begin{equation}\label{tensx}\begin{diagram}[nohug]TX\times TY &\rTo ^{\otimes}& T(X\times Y)\\
\uTo^{TX\times \eta_{Y}}&\ruTo_{t'_{X,Y}}&\\
TX \imes Y.&&&\end{diagram}\end{equation}
Now $t'_{X,Y}$ is itself 1-linear; but this 1-linearity may not be 
inherited by its 2-linear extension $\otimes$; it will be so if the 
monad is commutative, in the sense explained in Section \ref{ad10x}.

Similarly, $t''_{X,Y}:X\imes T(Y) \to T(X\times Y)$ extends uniquely 
over $\eta_{X}\times TY$ to a 1-linear $\tilde{\otimes}:T(X)\times 
T(Y) \to T(X\times Y)$, but its 2-linearity may not be inherited by 
its 1-linear extension.

An alternative description of the 2-linear $\otimes$ follows by 
specializing (\ref{2-lin-ext}); it is
$$\begin{diagram}TX\imes TY &\rTo ^{t''_{TX,Y}}&T(TX\imes 
Y)&\rTo^{T(t'_{X,Y})}&T^{2}(X\times Y)&\rTo^{\mu_{X\times 
Y}}&T(X\times Y)
\end{diagram}.$$
This is the description in \cite{MSMCC}, where $\otimes$ is called 
$\tilde{\psi}$ and $\tilde{\otimes}$ is called $\psi$.
The notion of commutativity of a strong monad was introduced in 
loc.cit.\ by the condition $\psi = \tilde{\psi}$.
Similarly, the 1-linear $\tilde{\otimes}$ is the composite
$$\begin{diagram}TX\imes TY &\rTo ^{t'_{X,TY}}&T(X\imes 
TY)&\rTo^{T(t''_{X,Y})}&T^{2}(X\times Y)&\rTo^{\mu_{X\times 
Y}}&T(X\times Y)
\end{diagram}.$$

Both $\otimes$ and $\tilde{\otimes}$ make the endofunctor $T$ into a 
monoidal functor. 
The reason why we call these two maps the ``Fubini'' maps is that 
when we below interpret ``elements'' in $T(X)$ (resp.\ in $T(Y)$) as 
(Schwartz) distributions, or as (Radon) measures, then the equality of them is a  form 
of Fubini's Theorem.

\section{The ``integration'' pairing}\label{ad7x} The following structure is 
the basis for the interpretation of $T(X)$ as a space of 
distributions or measures, in the ``double dualization'' sense of functional analysis.

Formally, the pairing consists in maps $T(X)\times (X\p B) \to B$, available 
for any $X\in \E$, and any $T$-algebra $B=(B,\beta )$. We denote it 
by a ``pairing'' bracket $\langle-,-\rangle$. Recall that we have the 
evaluation map $ev: X\imes (X\p B)\to B$, counit of the adjunction 
$(X\times - )\dashv ( X\p -)$. Then the pairing is the 1-linear extension 
of it over $\eta_{X }\times B$, thus $\langle-,-\rangle$ is 1-linear 
and makes the following triangle commute:
$$\begin{diagram}[nohug]T(X)\times (X\p B)&\rTo ^{\langle-,-\rangle}&B\\
\uTo^{\eta_{X}\times B}&\ruTo _{ev}&&\\
X\times (X\p B).
\end{diagram}$$

If $\phi: X\to B$ is an actual map in $\E$, it may be identified with 
a ``global'' element ${\bf 1}\to X\p B$, and we may consider the 
map $T(X)\to B$ given by $\langle -,\phi \rangle$. This map is then 
just the $T$-linear extension over $\eta_{X}$ of $\phi :X\to B$.

The brackets ought to be decorated with symbols $X$ and $(B, \beta 
)$, and it will quite evidently be natural in $B\in \E^{T}$ and in 
$X\in \E$, the latter naturality in the ``extranatural'' sense, 
\cite{Mac} 9.4, which we shall recall. If we use elements to express 
equations (a well known and rigorous technique, cf.\ e.g.\ \cite{SDG} 
II.1, even though objects in 
a category often are said to have no elements), the extraordinary 
naturality in $X$ is expressed: for any $f:Y \to X$, $\phi \in X\p B$ 
and $P\in T(Y)$,
\begin{equation*} \langle T(f)(P),\phi \rangle = 
\langle P, \phi \circ f \rangle ,
\end{equation*}
or more compactly, using $f_{*}$ to denote $T(f):T(Y)\to T(X)$ and 
$f^{*}$ to denote $f\p B$, i.e.\ precomposition with $f$,
\begin{equation}\label{extra-x}\langle f_{*}(P), \phi \rangle = 
\langle P, f^{*}(\phi )\rangle .
\end{equation}
To prove this equation, we observe that both sides are 1-linear in 
$P$, and therefore it suffices (Proposition \ref{one-x}) to prove that their precomposition 
with $\eta_{Y}\times B$ are equal; this in turn follows from 
naturality of $\eta$ w.r.to $f$, and the extraordinary naturality of 
$ev$.

\medskip
An alternative notation for the pairing:
$$\langle P, \phi \rangle = \int_{X} \phi \; dP = \int_{X} \phi(x)\; 
dP (x),$$
is motivated by the interpretation of $P\in T(X)$ as a measure (or as 
a  Schwartz distribution) on $X$, and of $\langle P, \phi \rangle$ as 
the integral of the function $\phi :X\to B$ w.r.to the measure $P$, 
or, as the value of the distribution on the ``test function'' $\phi$. 
In this context, the main case is where $B$ is the ``space of 
scalars'', but any other ``vector space'' (= $T$-algebra) $B$ may
be the value space for test functions. 
\section{The ``semantics'' map}\label{ad8x} By this, we understand the exponential transpose of the 
pairing map $T(X)\times (X\p B) \to B$; is is thus a map
$$\begin{diagram}T(X)&\rTo ^{\tau}& (X\p B)\p B\end{diagram}.$$
In this context, the $T$-algebra $B$ may be called ``the dualizing 
object''.
We could (should) decorate $\tau$ with symbols $X$ and/or $B$, if 
needed.
It is natural in $X\in \E$ (by extra-naturality of the pairing in the 
variable $X$), and it is natural in the dualizing object $B\in \E^{T}$, by the naturality 
of the pairing in $B\in \E ^{T}$. Why does it deserve the name 
``semantics''? 

Let us for a moment speak about $\E$ as if it were the category of 
sets. Thus, an element $P\in T(X)$ gets by $\tau$ 
interpreted as an $X$-ary operation that makes sense for any 
$T$-algebra $B$, and is preserved by morphisms of $T$-algebras $B\to 
B'$ in $\E ^{T}$. Here, we are talking about elements in $X\p B$ as 
``$X$-tuples of elements in $B$'', and functions $X\p B \to B$ are thus 
``$X$-ary operations'' that, given an $X$-tuple in $B$, returns a single 
element in $B$. Thus, $\tau_{B} (P)$ is the interpretation of the 
``syntactic'' (formal) $X$-ary operation $P$ as an actual $X$-ary operation on 
the underlying ``set'' of $B$; or, $\tau_{B}(P)$ is the {\em semantics} of $P$ in $B$, 
and it is preserved by any morphism of 
$T$-algebras. (This connection between monads on the category of 
sets, and universal algebra, - or equivalently, the idea of viewing 
monads on the category of sets as infinitary algebraic theories (in 
the sense of Lawvere) - goes back to the early days of monad theory 
in the mid 60's (Linton, Manes, Beck, Lawvere, \ldots). Our theory 
demonstrates that this connection makes sense for strong monads on 
cartesian closed categories.

\medskip

Part of the aim of the present study has been to demonstrate that the 
notion of an extensive quantity $P$ {\em distributed} over a given 
space $X$ can be encoded in terms of a strong monad $T$, but that it need not 
 be encoded in the ``double dualization'' paradigm of 
functional analysis, like the double dual $(X\p B)\p B$. In slogan 
form, ``{\em liberate distributions from the yoke of double 
dualization!}''. As a slogan, we shall partially ourselves contradict it, by 
arguing that ``there are enough test functions, provided we let $B$ vary'' 
(cf. Section \ref{ad12x} below).  On the other hand, one may push the slogan too far, if 
one reformulates it into ``{\em liberate syntax from the yoke of 
semantics!}'', which is only occasionally a fruitful endeavour. 
For, we hold that many distributions, say a 
distribution of a mass over a space (e.g.\  the distribution of tomato on a 
pizza) is not a syntactic entity, and does not in 
itself depend on the concept of double dualization. Mass 
distributions will be considered in Section \ref{ad15x} below.

\medskip 

\noindent{\bf Remark.} The $\tau$ considered here (with $B$ fixed) is actually a 
morphism of monads: even without an algebra structure on $B\in \E$, the functor 
$X \mapsto (X\p B)\p B$ is a strong monad on $\E$; if $T$ is another 
strong monad on $T$, we have the result that: to give $B$ a 
$T$-algebra structure is equivalent to give a strong monad morphism $\tau: T 
\Rightarrow (-\p B)\p B$, cf.\ \cite{DD}. The unit of the 
double-dualization monad
is the exponential transpose of the evaluation map $X\times (X\p 
B)\to B$. It is 
the ``Dirac delta map'' $\delta_{X} : X\to (X\p B)\p B$. In the category 
of sets: for  $x\in X$, $\delta_{X}(x)$ is the map $(X\p B) \to B$ which 
consists in evaluating at $x$.

 Since $\tau$ is the exponential 
transpose of the pairing, and the pairing was defined as 1-linear 
extension of evaluation, it follows that $\tau$ is the linear 
extension of the map $\delta$, in formula 
\begin{equation}\tau_{X }\circ \eta _{X} = \delta_{X},
\label{tau-delta}\end{equation}
or briefly, $\tau \circ \eta  = \delta$, ($B$ being a fixed $T$-algebra 
here). In \cite{CEQ}, we even denoted $\tau_{X}(x)$ by $\delta_{x}$.

\section{The object of $T$-linear maps}\label{ad9x}The following construction 
goes back to \cite{Bunge}. Let $(B, \beta)$ and $(C,\gamma)$ be two 
$T$-algebras, with $T$ a strong monad on $\E$. We assume that $\E$ 
has equalizers. Then we can out of the object $B\p C$ carve a 
subobject $B\p _{T}C$ ``consisting of'' the $T$-linear maps; 
precisely, it is the equalizer
$$\begin{diagram}
B\p _{T}C &\rTo & B\p C&\pile{\rTo ^{\beta \p C} \\ 
\rTo_{st. (TB\p \gamma )}}&T(B)\p C
\end{diagram}$$
in the lower map on the right, we compose forwards , so it is, in 
full
$$\begin{diagram}B\p C&\rTo ^{st_{B,C}}&T(B)\p T(C)&\rTo^{T(B)\p 
\gamma}&T(B)\p C,
\end{diagram}$$
expressing internally that a map $f:B\to C$ is a morphism of algebras if two 
particular maps $T(B)\to C$ are equal.

Now there is an evident notion of when a subobject $A'\subseteq A$ of (the 
underlying object of) an algebra $(A, \alpha )$ is a subalgebra. The 
object $B\p C$ inherits an algebra structure from that of $C$, using 
the cotensorial strength $\lambda$, so one may ask, is $B\p _{T}C$ a 
subalgebra ? This is not in general the case. If $B$ and $C$ are 
groups, the set of group homomorphisms $B\to C$ is not a subgroup of 
the group of all maps $B\to C$; it is so, however, for {\em commutative} 
groups. This leads to the topic of the next Section, commutativity of 
a monad.

\section{Commutative monads}\label{ad10x}

The notion of {\em commutative} monad (a strong monad with a certain {\em 
property}) was introduced in \cite{MSMCC}: it is a strong monad for 
which the two ``Fubini maps'' $\otimes$ and $\tilde{\otimes}$ agree. 
There are seve\-ral equivalent conditions describing commutativity, and 
we summarize these in the Theorem at the end of the section.

\begin{prop}\label{fivex}Let $T$ be commutative. Let $B=(B,\beta )$ and $C=(C,\gamma )$ be $T$-algebras, 
and assume that $f :X\times B \to C$ is 2-linear. Then its 1-linear 
extension $\overline{f}: T(X)\times  B \to C$ over $\eta _{X} \imes 
B$ is 2-linear (hence bilinear). (Similarly for 2-linear extensions of 1-linear maps.)
\end{prop}
{\bf Proof.} To prove 2-linearity of $\overline{f}$ means to 
prove commutativity of the following diagram (where the bottom line is 
$\overline{f}$, according to (\ref{1-lin-ext}))
$$\begin{diagram}TX \times TB&\rTo ^{t''}&T(TX\times B)&\rTo 
^{T(t')}&T^{2}(X\times B)&\rTo ^{T^{2}f}& T^{2}C&\rTo^{T\gamma}&TC \\
\dTo^{TX\times \beta }&&&&&&&&\dTo_{\gamma}\\
TX\times B &&\rTo_{t'}&T(X\times B)&\rTo_{Tf}&TC&\rTo_{\gamma }&&C.
\end{diagram}$$
The assumed 2-linearity of $f$ is expressed by the equation 
(composing from left to right)
\begin{equation}\label{2-lin-f}t''.Tf.\gamma = (X\times \beta ).f
\end{equation}
and the assumed commutativity of the monad is expressed by the 
equation
\begin{equation}\label{commxx} t''.Tt'.\mu _{X\times B}=t'.Tt''.\mu 
_{X\times B},
\end{equation}
using the explicit formulae for   $\otimes$ and $\tilde{\otimes}$ at 
the end of Section \ref{ad6x}.
We now calculate, beginning with the clockwise composite in the 
diagram:
$$\begin{array}{lcll}t''.Tt'.T^{2}f.T\gamma.\gamma&=&t''.Tt'.T^{2}f.\mu_{C}.\gamma&\mbox{ by an algebra law}\\
&=&t''.Tt'.\mu _{X\times B}.Tf.\gamma&\mbox{ by naturality of $\mu$}\\
&=&t'.Tt''.\mu_{X\imes B}.Tf.\gamma&\mbox{ by 
(\ref{commxx}), commutativity of the monad}\\
&=&t'.Tt''.T^{2}f.\mu_{C}.\gamma&\mbox{ by naturality of $\mu$}\\
&=&t'.Tt''.T^{2}f.T\gamma .\gamma&\mbox{ by an algebra law}\\
&=&t'.T(X\times \beta ).Tf.\gamma&\mbox{ by (\ref{2-lin-f}), 
2-linearity of $f$}\\
&=&(TX\imes \beta ).t'.Tf.\gamma&\mbox{ by naturality of $t'$,}
\end{array}$$
which is the counterclockwise composite of the diagram. This proves 
the Proposition.

\medskip
From  Proposition \ref{fivex} now immediately follows that when the 
monad is commutative, then $\otimes 
:T(X)\times T(Y) \to T(X\imes Y)$ is bilinear; for, it was 
constructed as the 2-linear extension of the 1-linear $t':T(X)\times 
Y \to T(X\imes Y)$.

It  follows rather easily that, conversely, bilinearity of $\otimes$ implies 
commutativity of the monad, cf.\ \cite{BCCM} Prop.\ 1.5. Thus, also, 
the extension properties in Proposition \ref{fivex} imply 
commutativity.

From Proposition \ref{fivex} also follows that when the 
monad is commutative, the pairing map
$T(X)\times (X\p B) \to B$ is bilinear; for, it was constructed as 
the 1-linear extension of the evaluation map $X \times (X\p B) \to 
B$, which is 2-linear. By passing to the 
exponential transpose,  this 2-linearity implies that its 
transpose, i.e.\ the semantics map $\tau :T(X) \to (X\p B)\p B$, factors 
through the subobject $(X\p B)\p _{T}B \subseteq (X\p B)\p B$. 

For commutative $T$, and for $B$ and $C$ arbitrary $T$-algebras, the 
subobject $B\p _{T }C \subseteq B\p C$ is actually a subalgebra (and 
in fact provides $\E^{T}$ with structure of a closed category), cf.\ 
\cite{CCGBCM}. 

The two last implications argued for here go the other way as well: 
bilinearity of the pairing implies commutativity; $B\p _{T}C$ a 
subalgebra of $B\p C$ (for all algebras $B$, $C$) implies commutativity. The 
proofs are not difficult, and will be omitted.

Finally, from \cite{BCCM}, Proposition 1.5 (v), we have the 
equivalence of commutativity with ``$\mu$ is monoidal''; also, 
commutativity implies that the monoidal structure $\otimes$ on the 
endofunctor $T$ is symmetric monoidal, cf.\ \cite{MSMCC} Theorem 
3.2).  We summarize:

\begin{thm}\label{cxx}Let $T=(T,\eta ,\mu )$ be a strong monad on $\E$. Then 
t.f.a.e., and define the notion of commutativity of $T$:

1) Fubini's Theorem  holds, i.e.\ the maps $\otimes$ and $\tilde{\otimes}: 
T(X)\times T(Y) \to T(X\times Y)$ agree, for all $X$ and $Y$ in $\E$.

2)  The map $\otimes :T(X)\times T(Y)\to T(X\times Y)$ is $T$-bilinear for all $X$ and $Y$.

3) For $(B,\beta ) $ and $(C, \gamma)$ algebras for $T$, and $X$ 
arbitrary, the 1-linear extension $T(X)\times B \to C$ of a 2-linear 
map $X \times B \to C$ is 2-linear (hence bilinear). (Similarly for 
2-linear extensions of 1-linear maps.)

4) For $(B,\beta ) $ and $(C, \gamma)$ algebras for $T$, the subspace 
$B\p _{T}C \subseteq B\p C$ is a sub-algebra.

5) For $(B,\beta )$ an algebra for $T$, and $X$ arbitrary, 
the ``semantics'' map $\tau : T(X) \to (X\p B)\p B$ factors through 
the subspace $(X\p B)\p _{T}B$.

6) For $(B,\beta )$ an algebra for $T$, and $X$ arbitrary, the
 pairing map $T(X)\times (X\p B) \to B$ is bilinear.

7) The monoidal structure $\otimes$ on the endofunctor $T$ is symmetric monoidal, 
and with respect to 
this monoidal structure $\otimes$, the natural transformation $\mu 
:T\circ T \Rightarrow T$ is 
 monoidal (hence the monad is a symmetric monoidal monad).
\end{thm}

\medskip

Recall the universal property of $X\times \eta_{Y}:X\times Y \to 
X\times T(Y)$ and of $\eta_{X}\times Y : X\times Y \to T(X)\times Y$ 
(Proposition \ref{one-x}). We have an analogous property
for $\eta_{X}\times \eta_{Y}:X\times Y \to T(X)\times T(Y)$: 
\begin{prop}\label{biextx}Let $T$ be a commutative monad. Let $B=(B, \beta )$ be a 
$T$-algebra. Then any $f:X\times Y \to B$ extends uniquely over $\eta 
_{X} \times \eta _{Y}$ to a bilinear $\overline{f}:T(X)\times T(Y) 
\to B$.
\end{prop}
{\bf Proof.} The extension is performed in two stages, along 
$\eta_{X}\times Y$, and then along $T(X)\times \eta _{Y}$. The first 
extension is unique, as a 1-linear map, the second is unique as a 
2-linear map. However, the 1-linearity of the first extension is 
preserved by the second extension, using clause 3) in the Theorem.

\medskip

Recall from Proposition \ref{onexx} that $t'_{X,Y} :T(X)\times Y \to 
T(X\imes Y)$ is an initial 1-linear map into  $T$-algebras, and 
similarly $t''_{X,Y}:X\times T(Y)\to T(X\imes Y)$ is an initial 
2-linear map. These properties join hands when $T$ is commutative:

\begin{prop}Let $T$ be a commutative monad. Then $\otimes :T(X)\times 
T(Y) \to T(X\times Y)$ is initial among bilinear maps to $T$-algebras.
\end{prop}
(Thus, $T(X\times Y)$ may be denoted $T(X)\otimes T(Y)$. Tensor 
products $A\otimes B$ for general $T$-algebras $A$ and $B$ 
may not exist, this depends on sufficiently 
many good coequalizers in $\E ^{T}$.)

\medskip
\noindent{\bf Proof.} Let $B=(B, \beta )$ be a $T$-algebra, and let 
$f: T(X)\times T(Y) \to B$ be bilinear.
 We must prove 
unique existence of a linear $\overline{f}$ making the right hand 
triangle in the following diagram commutative:
$$\begin{diagram}[nohug]X\imes Y&\rTo^{\eta_{X\times Y}}&T(X\times Y)\\
\dTo^{\eta_{X}\times \eta_{Y}}&\ruTo_{\otimes}&\dTo_{\overline{f}}\\
TX\imes TY&\rTo_{f}&B.
\end{diagram}$$
Since the upper triangle commutes, we see that $\overline{f}$ 
necessarily is the (unique) linear extension over $\eta_{X\times Y}$
of $(\eta_{X}\times \eta_{Y}).f$ (composing forwards).  To prove that this 
$\overline{f}$ does indeed solve the posed problem: The two maps 
$\otimes .\overline{f}$ and $f$ are 2-linear $TX\times TY \to 
B$, so by the universal property of $TX \imes \eta_{Y}$, to see that 
they agree, it suffices 
to see that the the agree when precomposed with $TX 
\times \eta_{Y}$.
Both of these two maps are  1-linear maps $TX \times Y \to B$, so by 
the universal property of $\eta _{X}\times Y$, it suffices to see 
that the two maps $X\imes Y \to B$ which one gets by precomposition 
with $\eta_{X}\times Y$ are equal. But they are both equal to 
$(\eta_{X}\times \eta_{Y}).f$, by construction of $\overline{f}$.

\section{Convolution; the scalars $R:=T({\bf 1})$}\label{ad11x} 
We henceforth assume that the monad $T$ is commutative.
If $m:X\times Y \to 
Z$ is a map in $\E$, we can use $\otimes$ to manufacture a map 
$T(X)\times T(Y) \to T(Z)$, namely the composite
$$\begin{diagram}TX \imes TY & \rTo ^{\otimes}&T(X\times Y)&\rTo 
^{T(m)}&TZ,
\end{diagram}$$  called the {\em convolution} along $m$. It is 
bilinear, since $\otimes$ is so, and $T(m)$ is linear. 
If $m:Y\times Y \to Y$ is associative, 
it follows from general equations for monoidal functors that the 
convolution along $m$ is associative. If $e: {\bf 1}\to Y$ is a two sided unit 
for the semigroup $Y,m$, then it gives rise to a unit for the 
convolution semigroup $T(Y)$, which thus acquires a monoid structure. If $m$ 
is furthermore commutative, then so is the convolution monoid; this 
latter fact depends on $T$ being a {\em symmetric} monoidal functor, which it is 
in our case, cf.\ clause 7) in the Theorem.

Similarly, if a monoid $(Y,m,e)$ acts on an object $X$ in an 
associative and unitary way, then the convolution monoid $T(Y)$ acts 
in an associative and unitary way on $T(X)$.

Now ${\bf 1}$ carries a unique (and trivial) monoid structure, and 
this trivial monoid acts on any object $X$, by the projection map 
$pr: {\bf 1}\times X \to X$. It follows that $T({\bf 
1})$ carries a monoid structure, and that this monoid acts on any 
$T(X)$.  Let us for the moment call this action of the {\em convolution 
action} on $T(X)$.
It is thus the composite
$$T\begin{diagram}{\bf 1}\times TX & \rTo^{\otimes} & T({\bf 1}\times X) \cong 
TX,\end{diagram}$$
it is bilinear, unitary and associative. It is the 2-linear extension 
over $T{\bf 1}\times \eta_{X}$ 
of $t'_{{\bf 1}, X}$ (followed by the isomorphism $T({\bf 1}\times 
X) \cong T(X)$); this follows from (\ref{tensx}). There is a similar 
description of a right action of $T({\bf 1})$ on $T(X)$, but $T({\bf 
1})$ is commutative, and these 
two actions agree. 

There is also an action by $T({\bf 1})$ on any 
$T$-algebra $B$, derived from the pairing, namely using the canonical 
isomorphism $i:B\to ({\bf 1}\p B)$,
$$\begin{diagram}T{\bf 1}\times B &\rTo^{T{\bf 1}\times i}&T{\bf 
1}\times ({\bf 1}\p B) &\rTo^{\langle-,-\rangle}&B.\end{diagram}$$ Let us call 
this the {\em pairing} action. 
\begin{prop}For a $T$-algebra of the form $T(X)$, the convolution 
action of ${T({\bf 1})}$ agrees with the pairing action.
\end{prop}
{\bf Proof.} Consider the diagram (where all the vertical maps are 
isomorphisms)
$$\begin{diagram}T{\bf 1}\times TX&\rTo^{t'}&T({\bf 1}\times 
TX)&\rTo^{T(t'')}&T^{2}({\bf 1}\times X)&\rTo^{\mu}&T({\bf 1}\times X)\\
\dTo^{T{\bf 1}\times i}&&\dTo^{T({\bf 1}\times 
i)}&&\dTo_{T^{2}(pr})&&\dTo_{T(pr)}\\
T{\bf 1}\times ({\bf 1}\p TX)&\rTo_{t'}&T({\bf 1}\times ({\bf 1}\p 
TX))&\rTo_{T(ev)}&T^{2}X&\rTo_{\mu}&TX.
\end{diagram}$$
The top line is the explicit form of (the 1-linear version of) 
$\otimes$, so the clockwise composite is the convolution action.
The counterclockwise composite is the pairing action. The left hand 
and right hand squares commute by naturality. The middle square comes 
about by applying $T$ to the square
$$\begin{diagram}{\bf 1}\times TX&\rTo ^{t''}& T({\bf 1}\times X)\\
\dTo^{{\bf 1}\times i}&&\dTo_{T(pr)}\\
{\bf 1}\times ({\bf 1}\p TX)&\rTo_{ev}&T(X).
\end{diagram}$$
We claim that this square commutes. In fact, both composites are 
2-linear, so it suffices to see that we have commutativity when we 
precompose by ${\bf 1}\times \eta_{X}$. This commutativity follows by 
naturality of $\eta$, together with the fact that $i$ and $pr$ 
are the exponential transposes of each other.

We conclude that the total diagram commutes, proving the desired 
equality of the two actions on an algebra $T(X)$ (the convolution 
pairing has only been defined for such algebras).

\medskip
(We note that the proof does not depend on commutativity of $T$, 
provided we take the 1-linear version of $\otimes$, as we did. Also, 
it is valid in any symmetric monoidal closed category, replacing 
${\bf 1}$ by the unit object, and $pr$ by the primitively given unit 
isomorphism $I\otimes X \to X$.)

\medskip

Since pairing is natural in $B\in \E ^{T}$, it follows that the pairing action
 has the property that any morphism in $\E^{T}$, i.e.\ any $T$-linear 
map $B \to C$, is equivariant. The Proposition then has the following 
Corollary:
\begin{prop}\label{equivarx}Any $T$-linear map $T(X)\to T(Y)$ is equivariant for the 
convolution action by $T({\bf 1})$.
\end{prop} 
The convolution action on the free $T$-algebra $T(X)$ by the monoid 
$T({\bf 1})$ is an associative 
unitary action, because $\otimes$ makes $T$ into a monoidal functor, as we 
remarked; one can probably prove that the pairing action by 
$T({\bf 1})$,  defined for all $T$-algebras, is likewise unitary and 
associative.

\medskip

Since the monoid $T({\bf 1})$ acts in a unitary, associative, and 
bilinear way on any $T(X)$, and $T$-linear maps preserve the action, 
it is reasonable to think of $T({\bf 1})$ as the (multiplicative) monoid of {\em 
scalars}, and to give it a special notation. We denote it by $R$;
$$R:=T({\bf 1});$$
 we 
denote its multiplication by a dot, and its unit by $1$. It is 
the map $\eta_{{\bf 1}}:{\bf 1}\to T({\bf 1})$.

(An addition on $R$ will be considered in Section \ref{ad16x} below.)

\subsection*{On totals}
For any $X\in E$, we have a unique map $!_{X}:X\to {\bf 1}$ to the 
terminal object. Composing it with the map $1:{\bf 1}\to R =T({\bf 
1})$ (i.e.\ by $\eta_{{\bf 1}}$) gives a map which we denote $1_{X}:X\to R$, the ``map with 
constant value 1''. 

Just by the fact that $T$ is a covariant functor, we have a map $\tot 
_{X}:T(X)\to R$, namely $T(!_{X}):T(X)\to T({\bf 1})=R$, ``formation 
of total''. From uniqueness of maps into the 
terminal, we have for any $f:X\to Y$  that 
$\tot_{X}= \tot_{Y}\circ T(f)$, or in elementwise terms, writing 
$f_{*}$ for $T(f)$:
$$\tot (P)=\tot (f_{*}(P)).$$
We have
\begin{equation}\label{totyy}\tot_{X}\circ \eta_{X} = T(!_{X})\circ \eta_{X} = \eta_{{\bf 
1}}\circ !_{X}\end{equation}
by definition of $\tot$ and by naturality of $\eta$; but the latter map is $1_{X}$.
-- In particular, we have, in elementwise terms, for $x\in X$,
\begin{equation}\label{totdelx}\tot (\eta_{X}(x))=1.\end{equation}


\begin{prop}\label{tottxx}The map $\tot _{X}:T(X)\to R$ equals the map $\langle 
-,1_{X}\rangle $.
\end{prop}

\noindent Expressed in terms of elements, for $P\in T(X)$,
$$\tot_{X}(P) = \langle P, 1_{X}\rangle = \int_{X} 1_{X}\; dP \in R$$
(recalling the alternative ``integral'' notation for the pairing).

\medskip
\noindent{\bf Proof.} The two maps $T(X)\to T({\bf 1})$ to be 
compared are $T$-linear, so it suffices 
to prove that their precomposites with $\eta_{X}$ agree. For $\tot 
_{X}$, (\ref{totyy}) shows that we get $1_{X}$.  For $\langle 
-,1_{X})$, the precomposition again gives  $1_{X}$ (cf.\ the description 
in Section \ref{ad7x}, immediately after the construction of the bracket).
\medskip

A consequence of the naturality of $\otimes$ w.r.to the maps $!:X\to 
{\bf 1}$ and $!:Y\to {\bf 1}$  is that \begin{equation}\label{tottensx}
\tot (P\otimes Q)=\tot (P)\cdot \tot (Q),
\end{equation}
where the dot denotes the product in $R=T({\bf 1})$.
Hence also, the convolution $P*Q$ of $P$ and $Q$ (along any map)
satisfies $\tot (P*Q)=\tot (P) \cdot \tot (Q)$.

\section{Schwartz distributions of compact support}\label{ad12x}
There exist cartesian closed categories $\E$ which contain the category of 
smooth manifolds, and also contain the category of convenient vector spaces 
(with smooth, not necessarily linear, maps), for instance the 
category $Lip^{\infty}$ of \cite{FK}. In the latter, the 
(convenient) vector space of Schwartz distributions of compact 
support on a manifold $X$ is represented as $(X\p {\mathbb R})\p 
_{T}{\mathbb R}$ for a suitable strong monad $T$ (with $T({\bf 1})= 
{\mathbb R})$, see \cite{FK} Theorem 5.1.1; this was 
one of the motivations for the present development. We shall not use 
this  material from functional analysis, except for terminology and 
motivation.
Thus, with $R= T({\bf 1})$,  $X\p R$ will be called  the {\em $T$-linear space of  (non-compact) 
test functions} on $X$, and $(X\p R)\p _{T}R$ is then the space of 
$T$-linear functionals on this space of test functions; $T$-linearity 
in the example unites  the two aspects, the algebraic (${\mathbb 
R}$-linear), and the topological/bornological. So in the \cite{FK} 
case of $Lip^{\infty}$, $(X\p 
R)\p _{T}R$ is  the {\em space of continuous linear functionals on 
the space of test functions}, which is how Schwartz 
distributions\footnote{often with ${\mathbb C}$ rather than ${\mathbb 
R}$ as dualizing object.} (of 
compact support) are defined, for $X$ a smooth manifold. 

Recall that for a  commutative monad $T$ on $\E$, we have the 
``semantics'' map $\tau _{X}:T(X) \to (X\p B) \p _{T}B$, for any 
$T$-algebra $B$, in particular, we have such a $\tau_{X}$ for $B=R=T({\bf 1})$,
\begin{equation}\label{scwartzx}T(X) \to (X\p R)\p _{T}R.\end{equation}
The functor $X\mapsto (X\p R)\p _{T}R$ is in fact itself a strong 
monad on $\E$, (not necessarily commutative), and $\tau$ is a strong 
natural transformation, and it is even a morphism of monads. (More 
generally, these assertions hold even when $R$ is replaced by any 
other $T$-algebra $B$.)

In the case of the $(\E$ , $T)$ of \cite{FK}, the map $\tau : (X\p R)\p 
_{T}R$ is an isomorphism for many $X$, in particular for smooth 
manifolds One may express this property verbally by saying that 
$T(X)$ is {\em reflexive} w.r.to $R$ in the category of $T$-algebras. 
For any $X$ in the $\E$ of \cite{FK}, it is a monic map; this one may express 
verbally: 
``there are enough ${\mathbb R}$-valued valued test functions to test equality of elements in 
$T(X)$''.

There are many  examples of $(\E ,T)$ where there are not enough 
$R=T({\bf 1})$-valued test functions ; thus there are 
interesting monads $T$ with $T({\bf 1}) ={\bf 1}$ (``affine monads'', 
see Section \ref{ad18x} below).

However,  we have the option of choosing other 
$T$-algebras $B$ as ``dualizing object''. Then there are enough test 
functions, in the following sense

\begin{prop}\label{testx}For any $X\in E$, there exists a $T$-algebra $B$ so that
$$\tau_{X}: T(X) \to (X\p B)\p _{T}B$$
is monic. The algebra $B$ may be chosen to be free, i.e.\ of the form 
$T(Y)$, in fact $T(X)$ suffices.
\end{prop} 
{\bf Proof.} Take $B:= T(X)$. We claim that the following diagram commutes: 
$$\begin{diagram}[nohug]T(X)&\rTo ^{\tau}&(X\p T(X))\p _{T}T(X)\\
&\rdTo_{=}&\dTo_{ev_{\eta}}\\
&&T(X)
\end{diagram}$$ 
where $ev_{\eta}$ denotes the map ``evaluate at the global element 
$\eta_{X}\in X\p T(X)$''. For, the maps to be compared are linear, so 
it suffices to see that they agree when precomposed with $\eta_{X}$.
Now $\eta _{X}$ followed by $\tau$ is the Dirac delta map $\delta 
:X\to (X\p B)\p B$ (cf.\ (\ref{tau-delta})), and 
 $\delta$ followed by ``evaluate at $\phi \in X\p B$'' is the map 
$\phi$ itself, by general ``$\lambda$-calculus''. So both maps give 
$\eta_{X}$ when precomposed by $\eta _{X}$.

\medskip  So not only are there enough test functions; for given $X$, 
even {\em one} single test function suffices, namely $\eta_{X}:X\to 
T(X)$.

\section{Multiplying a distribution by a function}\label{ad14x}
The action on distributions by (scalar valued) functions, to be 
described here, has as applications the notion of {\em density} of 
one distribution w.r.to another, and also the notion of {\em 
conditional} probability distribution.

Except for the subsection 
on ``totals", the considerations  so far are valid for any 
symmetric monoidal closed category $\E$, not just for a cartesian closed 
categories. The following, however, utilizes the universal property 
of the cartesian product.

Now that we are using that $\times$ is cartesian product, it becomes 
more expedient to use elementwise notation, since this will better allow us 
to ``repeat variables'',  like in the 
following. It is a description of a ``pointwise'' monoid structure on $X\p R$
derived from the monoid structure of $R$. The multiplication is thus 
a map $(X\p 
R)\times (X\p R) \to X\p R$ which in elementwise terms may be 
described, for $\phi$ and $\psi$ in $X\p R$,  and for $x\in R$,
$$(\phi \cdot \psi)(x):= \phi (x)\cdot \psi (x).$$ The unit is the 
previously described
$1_{X}$, in elementwise terms the map $x\mapsto 1 \in R$.
Similarly, the action of $R$ on $T(Y)$ for any $Y$ gives rise to a 
``pointwise'' action of $X\p R$ on $X\p T(Y)$.

We shall construct an action of $X\p R$ on the space $T(X)$. It has 
as a special case the ``multiplication of a (Schwartz-) distribution 
on $X$ 
by a scalar valued function $X\to R$'' known from classical 
distribution theory. We let the action be from the right, and denote 
it by $\vdash$:
$$\begin{diagram}T(X)\times (X\p R)& \rTo^{\vdash}& 
T(X),\end{diagram}$$
where $R= T({\bf 1})$; the $R$ here cannot be replaced by other 
$T$-algebras $B$.

To construct it, it suffices to construct a map $\rho :X\imes (X\p R) \to T(X)$ and 
extend it by 1-linearity over $\eta _{X}\times (X\p R)$. The map $\rho$ is constructed as the 
composite (with $pr_{1}$ being projection onto the first factor $X$ 
in the domain)
$$\begin{diagram}X\times (X\p T({\bf 1})) &\rTo^{\langle ev, pr_{1} 
\rangle}&T({\bf 1}) \times X &\rTo ^{t'_{{\bf 1},X}}& T({\bf 1}\times 
X)\cong T(X). 
\end{diagram}$$
(The map $\rho$ is 2-linear, and 2-linearity is preserved by the 
extension, by Theorem \ref{cxx}, so $\vdash$ is bilinear.) The 1-linear action $\vdash$ exists also without the 
commutativity assumption, but  then it cannot be asserted to be 
bilinear.)
The first map in this composite depends on the universal 
property of the cartesian product: a map into a product can be constructed 
by a pair of maps into the factors.

In elementwise notation, if $x\in X$ and $\phi \in X\p R$, we therefore have
\begin{equation}\label{dashetax}\eta (x)\vdash \phi = \phi (x)\cdot 
\eta (x)\end{equation}
(recalling the description (\ref{tensx}) of the action of $T({\bf 
1})$ on $T(X)$ in terms of $t'$).

In classical distribution theory, the following formula is the {\em 
definition} of how to multiply a distribution $P$ by a function 
$\phi$; in our context, it has to be {\em proved}. (The classical 
formula is the special case where $Y={\bf 1}$, so $T(Y)=R$.)

\begin{prop}\label{switchx}Let $P\in T(X)$, $\phi \in X\p R$ and $\psi \in  X\p 
T(Y)$. Then
$$\langle P\vdash \phi , \psi \rangle = \langle P, \phi \cdot \psi 
\rangle.$$ 
\end{prop}
{\bf Proof.} Since both sides of the claimed equation depend in a 
$T$-linear way on $P$, it suffices to prove the equation for the case 
where $P=\eta (x)$ for some $x\in X$. We calculate the left hand 
side, using (\ref{dashetax}):
$$\langle \eta (x)\vdash \phi ,\psi \rangle = \langle \phi (x)\cdot 
\eta (x), \psi \rangle = \phi (x) \cdot \langle \eta (x),\psi 
\rangle = \phi (x) \cdot \psi (x),$$
using for the middle equality that $T$-linearity implies equivariance w.r.to the action by 
$R$ (Proposition \ref{equivarx}). The right hand side of the desired 
equation similarly calculates
$$\langle \eta (x),\phi \cdot \psi \rangle = (\phi \cdot \psi 
)(x),$$
which is likewise $\phi (x)\cdot \psi (x)$, because of the pointwise character of the 
action of $X\p R$ on $X\p T(Y)$.

\medskip

\begin{corol}The pairing $\langle P, \phi \rangle$ (for $P\in T(X)$ 
and $\phi \in X\p R$) can be described in terms of $\vdash$ as 
follows:
$$\langle P, \phi \rangle = \tot (P\vdash \phi ).$$
\end{corol}
{\bf Proof.} Take $Y={\bf 1}$ (so $T(Y)=R$), and take $\psi = 1_{X}$. Then
$$\tot (P \vdash \phi)= \langle P \vdash \phi, 1_{X}\rangle = \langle 
P, \phi\cdot 1_{X}\rangle $$
using 
Proposition \ref{tottxx}, 
and then the Proposition \ref{switchx}. But $\phi\cdot 1_{X} 
= \phi$.

\medskip
Combining the ``switch'' formula in Proposition 
\ref{switchx} with Proposition \ref{testx} (``enough test 
functions''), we can derive properties of the action
 $\vdash$:
\begin{prop}The action $\vdash$ is associative and unitary. 
\end{prop}
{\bf Proof.} Let $\phi_{1}$ and $\phi_{2}$ be in $X\p R$ and let $P\in T(X)$. 
To see that $(P\vdash \phi _{1} )\vdash \phi_{2} = P\vdash (\phi_{1} 
\cdot \phi_{2} 
)$, it suffices by Proposition \ref{testx} to see that for any  
free $T$-algebra $B$ and any $\psi \in X\p B$, we have
$$\langle (P\vdash \phi _{1} )\vdash \phi_{2},\psi \rangle = 
\langle  P\vdash (\phi_{1} 
\cdot \phi_{2}),\psi \rangle,$$
but this is immediate using Proposition \ref{switchx} three times, and 
the associative law for the action of the monoid $X\p R$ on $X\p B$.
The unitary law is proved in a similar way.

\medskip 

To state the following result, we use for simplicity the notation 
$f_{*}$ for $T(f)$ (where $f:Y\to X$) and $f^{*}$ for $f\p B$, as in 
(\ref {extra-x}).
\begin{prop}[Frobenius reciprocity] For $f:Y\to X$, $P\in T(Y)$, and $\phi \in X\p 
R$,
$$f_{*}(P)\vdash \phi = f_{*}(P\vdash f^{*}(\phi )).$$
\end{prop}
{\bf Proof.} Both sides depend in a $T$-linear way of $P$, so it 
suffices to prove it for $P$ of the form $\eta_{Y}(y)$, for $y\in Y$.
We have
$$\begin{array}{lcll}f_{*}(\eta_{Y}(y))\vdash 
\phi&=&\eta_{X}(f(y))\vdash \phi&\mbox{\quad by naturality of $\eta$}\\
&=&\phi (f(y)) \cdot \eta_{X}(f(y))&\mbox{\quad by (\ref{dashetax})}\\
&=&\phi (f(y))\cdot f_{*}(\eta_{Y}(y))&\mbox{\quad by naturality of $\eta$}\\
&=&f_{*}\bigl( \phi (f(y))\cdot \eta_{Y}(y)\bigr)&\mbox{\quad by equivariance of 
$f_{*}$, Prop.\ \ref{equivarx}}\\
&=& f_{*}
\bigl( (\eta_{Y}(y))\vdash f^{*}(\phi )\bigr)&\mbox{\quad  by $f_{*}$ 
applied to 
(\ref{dashetax}).}
\end{array}$$
An alternative proof is by using the ``enough test-functions'' 
technique.

\subsection*{Density functions}

For $P$ and $Q$ in $T(X)$, it may happen that $Q= P\vdash \phi$ for 
some $\phi \in X\p R$, in which 
case one says that $Q$ has a {\em density} w.r.to $P$, namely the 
(scalar-valued) function $\phi$. Such $\phi$ may not exist, and  if it does, it may 
not be unique. In the case of (non-compact) Schwartz distributions, 
one case is particularly important, namely where $X$ is ${\mathbb 
R}^{n}$, and $P$ is Lebesgue measure; then if $Q$ has a density 
function $\phi$ w.r.to $P$, one sometimes identifies the distribution 
$Q$ with the function $\phi$. Such identification, as stressed by 
Lawvere, leads to loss of the distinction between the covariant 
character of $T(X)$ and the contravariant character of $X\p R$, more 
specifically, between extensive and intensive quantities.

\section{Mass distributions, and other extensive 
quantities}\label{ad15x}
An extensive quantity of a given type $m$ may, according to Lawvere, be 
modelled mathematically by a covariant functor $M$ from a category $\E$ of 
spaces to a ``linear'' or ``additive'' category $\A$, with suitable structure and properties. 
In particular, the category $\E$ should be ``lextensive''  
(cartesian closed categories with finite coproducts have this 
property), and $\A$ should be ``linear'' (categories of modules over 
a rig\footnote{commutative semiring with unit} have this property). In such a 
category, any object is canonically an abelian semigroup. Also, the 
functor $M$ should take finite coproducts to bi-products, by a 
certain canonical map.

It is a reasonable mathematical model that the quantity type of {\em mass} 
(in the sense of physics) is such a functor, with $\A$ being the 
category of modules over the rig of non-negative reals. If $M(X)$ 
denotes the set of possible mass distributions over the space $X$, 
then if $P_{1}$ and 
$P_{2}$ are masses
which are distributed over the  space $X$, one may, almost by physical 
construction,  combine them into one mass, distributed over $X$; so 
$M(X)$ acquires an additive structure; similarly, one may re-scale a 
mass distribution by  non-negative reals.
 The covariant functorality of $M$ has for its germ the idea that for 
a mass $P$ distributed over a space $X$, one can form its {\em 
total} $\tot (M)$. Such total mass may be something like 100 gram, so is not 
in itself a scalar $\in {\mathbb R}_{+}$, but only becomes so after 
choice of a {\em unit} mass, like gram, so the scalar is not {\em 
canonically} associated to the total mass.

{\em Intensive} quantities are derived from extensive ones as 
densities, or {\em ratios} 
between extensive ones. Thus, ``specific weight'' is an 
intensive quantity of type $m\cdot v^{-1}$ (where $v$ is the quantity 
type of volume (if the space $X$ is a block of rubber, its volume 
(Lebesgue measure)  may 
vary;  so volume is, like mass, a distribution varying over $X$).

Here, we shall only be concerned with {\em pure} intensive quantities; 
they are those of type $m\cdot m^{-1}$ for some extensive quantity type $m$. 
For the example of mass, an intensive quantity of type $m\cdot m^{-1}$ over the 
space $X$ may be identified with a function $X\to {\mathbb R}_{+}$. Other extensive quantity types may 
naturally give rise to other spaces of scalars; e.g.\ the quantity 
type of electric charge gives, in the simplest mathematical model, 
rise to the space ${\mathbb R}$ of {\em all} real numbers (since 
charge may also be negative), or, in a more sophisticated model, 
taking into account that charge is considered to be something 
quantized (so can be counted by integers), the space ${\mathbb Q}$ 
should in principle be the correct one.

So every quantity type defines in principle its own rig $R$ of scalars. 
However, it is a remarkable fact that so many give rise to ${\mathbb 
R}$, as the simplest mathematical model.

  In real life, pure 
quantities come after the extensive physical ones.  The theory presented in the 
preceding Sections  
follows, however, the mathematical tradition\footnote{There is a, 
minor, mathematical tradition in the other direction, namely 
Geometric Algebra, where on {\em constructs} a ring of scalars {\em from} 
 geometry.},  and puts the  theory of {\em pure} 
quantities as the basis. 

So to use Lawvere's concepts, our cartesian closed category $\E$ is 
the category of spaces; for $X$ a space, $T(X)$ is the space of 
 pure extensive quantities, varying over $X$; $R=T({\bf 1})$ is the space 
of scalars, and $X\p R$ is the space of intensive pure quantities 
varying over $X$. The functor $T$ may be seen as taking values in the 
category $\A:=\E ^{T}$ of $T$-algebras; it is not quite a linear 
category, but at least, all objects carry an action by the monoid of 
scalars, and all maps are equivariant with respect to these, cf.\ 
Proposition \ref{equivarx}. 
(Additive structure in $\A$ will be considered in the next Section.)

 \medskip

The aim of the remainder of the present Section is to demonstrate how 
e.g.\ the theory of a non-pure extensive quantity type, like mass, 
embeds into the theory considered so far; more precisely, to give a 
framework which vindicates the idea that ``as soon as a unit of mass 
is chosen, there is no difference between the theory of mass 
distributions, and the theory of distributed pure quantities''.
The framework is the following:

Let $T$ be a commutative monad on $\E$. Consider 
another strong endofunctor $M$ on $\E$, equipped with an action $\nu$ 
by $T$,
$$\nu : T(M(X))\to M(X)$$
strongly natural in $X$, and with $\nu$ satisfying a unitary and associative law. Then 
every $M(X)$ is a $T$-linear space by virtue of $\nu _{X}:T(M(X))\to 
M(X)$, and morphisms of the form $M(f)$ are $T$-linear.
Let $M$ and $M'$ be strong endofunctors equipped with such 
$T$-actions.
There is an evident notion of when a strong natural transformation 
$\lambda  :M\Rightarrow M'$ is compatible with the $T$-actions, so we 
have a category of $T$-actions. The endofunctor $T$ itself is an 
object in this category, by virtue of $\mu$. We say that $M$ is a 
{\em $T$-torsor} if it is isomorphic to $T$ in the category of 
$T$-actions. Note that no particular such isomorphism is chosen.

Our contention is that the category of $T$-torsors is a mathematical  model of 
(not necessarily pure) quantities $M$ of type $T$ (which is the 
corresponding pure quantity).

The following Proposition expresses that isomorphisms of actions $\lambda: T\cong M$ are 
determined  by $\lambda _{{\bf 1}}:T({\bf 1})\to M({\bf 1})$; in the example, the 
latter data means: choosing a {\em unit} of mass.

\begin{prop} If $g$ and $h:T\Rightarrow M$ are isomorphisms of $T$-actions, and if 
$g_{{\bf 1}}=h_{{\bf 1}}:T({\bf 1})\to M({\bf 1})$, then 
$g=h$.\end{prop}
{\bf Proof.} By replacing $h$ by its inverse $M\Rightarrow T$, it is clear 
that it suffices to prove that if $\rho :T\Rightarrow T$ is an isomorphism of 
$T$-actions, and $\rho_{{\bf 1}}= id_{T({\bf 1})}$, then $\rho$ is the identity 
transformation. As a morphism of $T$-actions, $\rho$ is in particular 
a {\em strong} natural transformation, which implies that the right hand  
square in the 
following diagram commutes for any $X\in \E$; the left hand square commutes 
by assumption on $\rho _{{\bf 1}}$:
$$\begin{diagram}X\times {\bf 1}&\rTo^{X\times \eta_{{\bf 1}}}&X\times 
T({\bf 1})&\rTo^{t''}&T(X\times {\bf 1})\\
\dTo^{=}&&\dTo^{X\times \rho_{{\bf 1}}}&&\dTo_{\rho_{X\times {\bf 1}}}\\
X\times {\bf 1}&\rTo_{X\imes \eta_{{\bf 1}}}&X\imes T({\bf 1})&\rTo _{t''}&T(X\imes {\bf 1})
\end{diagram}$$ 
Now both the horizontal composites are $\eta _{X\times {\bf 1}}$, by general 
theory of tensorial strengths. Also $\rho _{X\times {\bf 1}}$ is 
$T$-linear. Then uniqueness of $T$-linear extensions over 
$\eta_{X\times {\bf 1}}$ 
implies that the right hand vertical map is the identity map. Using the natural 
identification of $X\times 
{\bf 1}$ with $X$, we then also get that $\rho_{X}$ is the identity map of 
$T(X)$.

\section{Additive structure}\label{ad16x}For the applications of distribution 
theory, one needs not only that distributions on a space can be 
multiplied by scalars, but also that they can be {\em added}. In our 
context, this will follow if the category $T$-algebras is an {\em additive} 
(or {\em linear}) category, in a sense which we shall recall. Most textbooks, e.g.\ 
\cite{Mac}, \cite{Borceux} define the notion of additive category  
(with biproducts)  as a property of 
categories $\A$ which are already equipped with the {\em structure} of enrichment in 
the category of abelian monoids (or even abelian groups). However, we shall need that such 
enrichment structure derives canonically from a certain {\em property} of the 
category. This is old wisdom, e.g.\ described in Pareigis' 1969 
textbook \cite{Pareigis}, Section 4.1.

We recall briefly the needed properties: Consider a category $\A$ 
with finite products and finite coproducts. If the unique map from 
the initial object to the terminal object is an isomorphism, this 
object is a zero object $0$. The zero object allows one to define a 
canonical map $\delta$ from the coproduct  of $n$ objects $B_{i}$ to their 
product. If this $\delta$ is an isomorphism for all $n$-tuples of 
objects, this coproduct (and also, the product) is a biproduct 
$\oplus _{i}B_{i}$.
Then every object acquires the structure of an abelian monoid object 
in $\A$, with the codiagonal map $B\oplus B \to B$ as addition, and 
with the map unique map  $0\to B$ as unit. This in turn implies a 
canonical enrichment of $\A$ in the category of abelian monoids. We refer to \cite{Pareigis}, (except that 
loc.cit.\ item 4) includes a property which implies subtraction of 
morphisms, i.e.\ an enrichment in the category of abelian groups, but 
so long subtraction can be dispensed with, then so can item 4).)

\medskip

Consider now a category $\E$ with finite products and finite 
coproducts, and consider a monad $T$ on $\E$. Let $F: \E \to \E ^{T}$ 
be the functor $X\mapsto (T(X),\mu_{X})$, i.e.\ the functor 
associating to $X$ the free $T$-algebra on $X$. If $T(\emptyset )$ is 
the  terminal object 
${\bf 1}$ (where $\emptyset$ is the initial object), then it is clear that $F(\emptyset )$ is a zero object in 
$\E ^{T}$. Even though $\E^{T}$ may not have all finite coproducts, 
at least it has finite coproducts of free algebras, and so one may 
consider the canonical map $\delta :F(X)+F(Y) \to F(X)\times F(Y)$ (whose 
underlying map in $\E$ is a map $T(X+Y)\to T(X)\times T(Y)$,  
easy to describe directly, using $0$ and the universal property of 
$\eta_{X+Y}$). If it is an isomorphism for all $X$ and 
$Y$, the category of free algebras will therefore  be additive with 
biproducts $\oplus$. The addition map for an object $T(X)$ is thus 
the inverse of $\delta$ followed by $T(\nabla )$, where $\nabla 
:X+X\to X$ is the codiagonal. Since $\delta$ is $T$-linear, hence so 
is its inverse. Also $T(\nabla )$ is $T$-linear, so the addition map
$T(X)\times T(X) \to T(X)$ is $T$-linear.
For a more complete description, see \cite{CJ} 
or \cite{MEQ}. (In \cite{CJ}, it is proved that in fact that the 
category of {\em all} $T$-algebras, i.e.\ $\E^{T}$, is 
additive with biproducts.)

Thus, in particular, any $T$-linear map will be additive. It is 
tempting to conclude ``hence any bilinear (i.e.\ $T$-bilinear) map will be 
bi-additive''. This conclusion is true, and similarly 1- or 2-linear will be 
additive in the first (resp.\ in the second) variable; but in all 
three cases, an argument is required. Let us sketch 
the proof of this for the case of 2-linearity. Since $t''_{X,Y}$ is 
initial among 2-linear maps out of $X\times T(Y)$ (Proposition 
\ref{onexx}), it suffices to 
prove that $t''$ is additive in the second variable, which is to say 
that the following diagram commutes (where $\tilde{\Delta}$ is 
``diagonalizing $X$ and then taking middle four interchange''):
$$\begin{diagram}X\imes TY\times TY&\rTo^{\tilde{\Delta}}&X\imes 
TY\times X\times TY&\rTo ^{t''\times t''}&T(X\times Y)\times 
T(X\times Y)\\
\dTo^{X\times (+)}&&&&\dTo_{(+)}\\
X\imes TY&&\rTo_{t''}&&T(X\times Y).
\end{diagram}$$
To see this commutativity, it suffices to see that the two composites 
agree when precomposed with 
\begin{equation}\label{hintx}\begin{diagram}X\imes TY&\rTo^{X\times T(in_{i})}&X\times 
T(Y+Y)&\rTo^{\delta}_{\cong}&X\times TY\times TY,
\end{diagram}\end{equation}
($i=1,2$); it is easy to see that we get $t''$ in both cases (and for 
$i=1$ as well as for $i=2$). (Hint: use that (\ref{hintx}) 
postcomposed with $X\times (+)$ is the identity map on $X\times TY$, 
cf.\ the construction \cite{MEQ}, equation (32).) 

Let us summarize:

\begin{itemize}
\item the addition map $T(X)\times T(X) \to T(X)$ is $T$-linear
\item $T$-linear maps are additive
\item $T$-bilinear maps (resp.\ $T$-1-linear, resp.\ $T$-2-linear 
maps) are bi-additive (resp.\ additive in the first variable, resp. 
additive in the second variable)
\end{itemize} 

\noindent It follows that
the monoid of scalars $R=T({\bf 1})$ carries a $T$-linear addition 
$+:R\times R \to R$, and that the $T$-bilinear multiplication 
$R\times R \to R$ is bi-additive. So $R$, with its multiplication 
and addition, is a rig. Similarly, the multiplicative action of $R$ on 
any $T(X)$ is bi-additive. Since any $T$-linear map $T(X)\to T(Y)$ is 
equivariant for the action, and, as a $T$-linear map, is additive, it 
follows that any $T(X)$ is a module over the rig $R$, and that 
$T$-linear maps $T(X)\to T(Y)$ are $R$-linear, in the sense of module 
theory. 
\medskip

Since each object in $\E^{T}$ carries a canonical structure of 
abelian monoid, one may ask whether these abelian monoids have the 
{\em property} of being abelian {\em groups}. This is then a property 
of the monad. It is, more compactly, equivalent to the property: 
$1\in R$ has an additive inverse $-1$.
If this is so, difference-formation, for any $T$-algebra $A$, will be a 
$T$-linear map $A\times A \to A$.

We shall need this further property in Section \ref{ad19x}: 
differential calculus depends on having differences.


\section{Distributions on the line $R$  of scalars}\label{ad17x}
We now consider distributions $P\in T(X)$, where $X$ is the space $R$ of scalars.
Then $X\p R = R\p R$ has some particular (global) elements, to which 
we can apply $\langle P,-\rangle$, namely the monomials $x^{n}$, and 
also, $R $ has an addition $R\imes R \to R$, and we have convolution 
along this, which we shall denote $*:T(R)\times T(R) \to T(R)$.
 (There are some quite evident generalizations to the 
space $X=R^{n}$, but for simplicity, we stick to the 1-dimensional 
case.) Examples of elements in $T(R)$ are: probability distributions 
of random variables, cf.\ Section \ref{ad18x}.

The {\em $n$th moment}  of $P\in T(R)$ is defined as $\langle P, 
x^{n}\rangle \in R$, where $x^{n}:R\to R$ is the map $x\mapsto x^{n}$. In 
particular, the $0$th moment of $P$ is $\langle P, 1_{R} \rangle$, 
which we have considered already under the name $\tot (P)$, cf.\  
Proposition \ref{tottxx}. 
The $1$st moment is $\langle P, x \rangle$, i.e.\ $\langle P, 
id_{R}\rangle$, or, in the ``integral'' notation also used for the 
pairing,
$\int_{R}x\; dP(x)$. When $P$ is a probability distribution, this 
scalar is usually called the {\em expectation} of $P$. We therefore also 
denote it $E(P)$, so $E(P):=\langle P, x \rangle$, for arbitrary $P 
\in T(R)$.
(We shall not here consider $n$th moments for $n\geq 2$.)

 Recall that 
$R:=T({\bf 1})$, so $T(R) = T^{2}({\bf 1})$, and we therefore have 
the map $\mu_{{\bf 1}}:T(R)\to R$.
\begin{prop}\label{mu-exp}The ``expectation'' map $E:T(R)\to R$ equals $\mu_{{\bf 
1}}:T^{2}({\bf 1})\to T({\bf 1})$.
\end{prop}
{\bf Proof.} To prove $\langle P, id_{R}\rangle = \mu_{{\bf 1}}(P)$, 
we note that both sides of this equation depend in a $T$-linear way of 
$P$, so it suffices to prove it for the case where $P=\eta_{R}(x)$ 
for some $x\in R$. But $\langle \eta_{R}(x), id_{R}\rangle$ is 
``evaluation of $id_{R}$ on $x$'', by construction of the pairing, so 
gives $x$ back. On the other hand $\mu_{{\bf 1}}\circ \eta_{R}$ is 
the identity map on $R$, by a monad law, so this composite likewise 
returns $x$.

\medskip

Recall from (\ref{totyy}) that $\tot_{X} \circ \eta _{X}= 1_{X}$ as maps from 
$X$ to $R$ (for any $X$). From 
Proposition \ref{mu-exp}, we 
have that $E \circ \eta_{R} =id_{R}$; in elementwise terms, for $x\in 
R$,
\begin{equation}\label{expdelx}\tot (\eta (x))=1\quad \mbox{ ; } 
\quad E(\eta (x)) = x,\end{equation} 
(where $\eta$ denotes $\eta _{R}$). Using these two facts, we can 
prove
\begin{prop}\label{newx} For $P$ and $Q$ in $T(R)$, we have
$$E(P*Q)= E(P)\cdot \tot (Q) + \tot(P)\cdot E(Q).$$
\end{prop}
{\bf Proof.} Convolution is $T$-bilinear, and $E$ and $\tot$ are $T$-linear. Since 
addition is a $T$-linear map, it follows that both sides of the 
desired equation depend in a $T$-bilinear way on $P,Q$. Therefore, by 
Proposition \ref{biextx}, it 
suffices to prove it for the case where $P=\eta (x)$ and $Q=\eta(y)$ 
(with $x$ and $y$ in $R$).
Then $E(P)=x$, $E(Q)=y$ and $\tot (P)=\tot (Q)=1$. Also, since $*$ is 
convolution along $+$, we have
$E(P*Q)=E(\eta (x)*\eta (y)) = E(\eta (x+y)) = x+y$. Then  the 
result is immediate. 

\medskip

 Let $h$ be a homothety on $R$, $x\mapsto b\cdot x$, for some 
$b \in R$. It 
is a $T$-linear map, since the multiplication on $R$ is $T$-bilinear. 
We claim 
\begin{equation}E(h_{*}(P)) = h(E(P))\label{homothx}\end{equation}
 for any $P\in T(R)$. To see this, 
note that both sides depend in a $T$-linear way on $P$, so it 
suffices to see it for the case of $P=\eta (x)$ ($x\in R$). But
$E(h_{*}(\eta (x))) = E (\eta (h(x))) =h(x)=  h(E(\eta (x)))$.

We also consider the effect of a translation map $\alpha$ on $R$, 
i.e.\ a map of the form $x\mapsto x+a$ for some $a\in R$. We claim
\begin{equation}\label{translxx}
\alpha_{*}(P) = P*(\eta (a))\end{equation} for 
any $P\in T(R)$. To see this, we note that both sides 
here depend in a $T$-linear way on $P$, so it suffices to see it for 
$P=\eta (x)$. Then both sides give $\eta (x+a)$.

Applying Proposition \ref{newx} to $P*\eta (a)$, we therefore 
have, for $\alpha =$ translation by $a\in R$,
$$E(\alpha _{*}(P))= E(P) + \tot (P)\cdot a.$$
In particular, if $P$ has total 1,
\begin{equation}\label{translx}E(\alpha _{*}(P)) = E(P) + a = \alpha 
(E(P)).\end{equation}

Recall that an affine map $R\to R$ is a map $f$ which can be written 
as a homothety followed by a translation. Combining (\ref{homothx}) 
and (\ref{translx}), we therefore get
\begin{prop}\label{CG1x}Let $f:R\to R$ be an affine map. Then if $P\in T(R)$ has 
total 1, we have
$$E(f_{*}(P)) = f(E(P)).$$
\end{prop}
If $P \in T(R)$ has $\tot (P)$ multiplicatively invertible, we define 
{\em the center of gravity} $cg(P)$ by
$$cg(P):= E((\tot (P))^{-1}\cdot P).$$
Since $E$ is $T$-linear, it preserves multiplication by scalars, so 
$Q:=(\tot (P))^{-1}\cdot P$ has total 1, and the previous Proposition 
applies. So for an affine $f:R\to R$,
$E(f_{*}(Q)) = f(E(Q))$. The right hand side is $cg(P)$, by 
definition. But $f_{*}$ preserves formation of totals (any map does), 
and then it is easy to conclude that the left hand side is 
$cg(f_{*}(P))$. Thus we have
\begin{prop} Formation of center of gravity is preserved by affine 
maps $f$, $cg(f_{*}P)=f(cg(P))$.
\end{prop}
Thus, our theory is in concordance with the truth that ``center of gravity for a mass distribution on a line does not depend on choice of 0 and 
choice of unit of length''.

\section{The affine submonad; probability distributions}\label{ad18x}
A strong monad $T$ on a cartesian closed category $\E$ is called {\em 
affine} if $T({\bf 1})={\bf 1}$. 
For algebraic theories (monads on the category of sets), this was introduced in \cite{Wraith}. For 
strong  monads, it was proved in \cite{BCCM} that this is equivalent 
to the assertion that for all $X,Y$, the map
$\psi _{X,Y}:T(X)\times T(Y) \to T(X\imes Y)$ is split monic with 
$\langle T(pr_{1}),T(pr_{2})\rangle :T(X\times Y) \to T(X)\times 
T(Y)$ as retraction. 
In \cite{Lindner}, it was proved that if $\E$ has finite limits, any 
 commutative monad $T$ has a maximal affine submonad $T_{0}$, the 
``affine part of $T$''. It is likewise a commutative monad.
Speaking in 
elementwise terms, $T_{0}(X)$ consists of those 
distributions whose total is $1\in T({\bf 1})$. We consider in the 
following a commutative monad $T$ and its affine part $T_{0}$.

Probability distributions   have by definition total $1\in R$, and  
 take values in the interval from 0 to 1, in the sense that $0\leq \langle 
P, \phi \rangle \leq 1$ if $\phi$ is a multiplicative idempotent in 
$X\p R$. We do not in the present 
article consider any order relation on $R$, so there is no ``interval 
from 0 to 1''; so we are stretching 
terminology a  bit when we use the word ``probability distribution on $X$'' for the 
elements of $T_{0}(X)$, but we shall do so.

 If $P\in T_{0}(X)$ and $Q\in 
T_{0}(Y)$, then $P\otimes Q\in T_{0}(X\times Y)$, cf.\ 
(\ref{tottensx}) and the remark following it. (It
also can be seen from the fact that the inclusion of strong monads $T_{0}\subseteq T$ is 
compatible with the monoidal structure $\otimes$.)
From this in turn follows that e.g.\  probability 
distributions are stable under convolution.

Most of the theory developed presently works for the monad $T_{0}$ 
just as it does for $T$; however, since $T_{0}({\bf 1})={\bf 1}$ is 
trivial (but, in good cases, $T_{0}({\bf 2}) = T({\bf 1})
$), 
the theory of  ``multiplying a distribution by a 
function'' (Section \ref{ad14x}) trivializes. So we shall consider 
the $T_{0}$ in conjunction with $T$, and $R$ denotes $T({\bf 1})$.
Clearly, $T_{0}(X)$ is not stable under multiplication by scalars 
$\in R$; for, the total gets multiplied by the scalar as well. In 
particular, we cannot, for $P\in T_{0}(X)$, expect that $P\vdash \phi$ 
is in $T_{0}(X)$. Nevertheless $P\vdash \phi$ has a probabilistic 
signi\-ficance, provided $\langle P,\phi \rangle $ is multiplicatively 
invertible in $R$. If this is the case, we may form  a distribution denoted 
$P\! \mid \! \phi$,
$$P\mid \phi := \langle P,\phi \rangle ^{-1} \cdot (P\vdash \phi ) = 
\lambda \cdot (P\vdash \phi ),$$
(writing $\lambda 
\in R$ for $\langle P,\phi \rangle ^{-1}$) and so (using Proposition \ref{switchx}, and bilinearity) we have
$$\tot (P\! \mid \! \phi) =  \langle P\! \mid \! \phi , 1_{X}\rangle 
=\lambda \cdot \langle P\vdash \phi 
,1_{X}\rangle =\lambda \cdot \langle P, 
\phi \cdot 1_{X}\rangle  =\lambda \cdot \langle P, \phi \rangle =1,$$
so that $P\! \mid \! \phi$ again is a probability distribution.
More generally $P\! \mid \!  \phi$,  satisfies, for any $\psi \in X\p R$,
\begin{equation}\label{condix}\langle P\! \mid \! \phi ,\psi\rangle = \lambda \langle P, \phi \cdot 
\psi \rangle .\end{equation} Now,  multiplicatively idempotent elements in $X\p 
R$ may reasonably be seen as {\em events} in the {\em outcome space} 
$X$, and if $\phi$ and $\psi$ are such events, also $\phi \cdot \psi$ 
is an event, and it deserves the notation $\phi \cap \psi$ (simultaneous 
occurrence of the two events). In this 
notation, and with the defining equation for $\lambda$, the equation 
(\ref{condix}) reads
$$\langle P\! \mid \! \phi ,\psi \rangle = \frac{\langle P, \phi \cap \psi 
\rangle}{\langle P,\phi \rangle},$$
so that $P\! \mid \! \phi$ is ``the conditional probability that $\psi$ 
(and $\phi$) 
occur, given that $\phi$ does''.

\medskip
\subsection*{Random variables and their distributions}
A random variable  $X$ on an outcome space $\Omega$ defines a probability 
distribution $\in T_{0}(R)$, and is often identified with this 
distribution. A {\em pair} $X_{1},X_{2}$ of random variables on the same $\Omega$, on the other hand, 
cannot be identified with a pair of distributions on $R$, but rather 
with a distribution $P\in T_{0}(R^{2})$, the {\em joint} probability 
distribution of the two random variables. This $P$  {\em gives rise} to a 
pair of distributions $P_{i}:= (pr_{i})_{*}(P)$ ($i=1,2)$ (the {\em 
marginal distributions} of $P$); so we have 
$(P_{1},P_{2})\in T(R)\times T(R)$. We may compare
$P$ with $P_{1}\otimes P_{2}$. To say that $P= P_{1}\otimes P_{2}$ is 
expressed by saying that the two random variables $X_{1}$ and $X_{2}$ are {\em independent}.

(On the other hand, the marginal distributions of $P\otimes Q$, for 
$P$ and $Q \in T_{0}(R)$, agree with $P$ and $Q$; this is a general property 
of affine monads, mentioned earlier.)

A pair of random variables $X$ and $Y$ on $\Omega$ has a sum, which 
is a new random variable; it is often denoted $X+Y$, but the 
probability distribution $\in T(R)$ of this sum is not a sum of two 
probability distributions (such a sum would anyway have total 2, 
not total 1). If $X$ and $Y$ are independent, the sum $X+Y$ has as 
distribution the convolution along $+$ of the two individual 
distributions. In general, the distribution of $X+Y$ is $+_{*}(P)$ 
where $+: R^{2}\to R$ is the addition map and $P \in T(R^{2})$ is the 
joint distribution of the two random variables.

In text books on probability theory (see e.g.\ \cite{Cramer}), one sometimes sees the formula
$$E(X+Y)=E(X)+E(Y)$$
for the expectation of a sum of two random variables on $\Omega$ (not 
necessarily independent). The meaning of this non-trivial formula 
(which looks trivial, because of the notational confusion between a 
random variable and its distribution) is, 
in the present framework, a property of distributions on $R$. Namely 
it is a way to record the commutativity of the diagram
$$\begin{diagram}T(R\times R)&\rTo^{T(+)}&T(R)\\
\dTo^{\beta}&&\dTo_{E}\\
R\times R&\rTo_{+}&R.\end{diagram}$$
Here, $\beta$ is the $T$-algebra structure on $T(R^{2})$, i.e.\ 
$\mu_{2}:T^{2}(2) \to T(2)$ (recall that $2={\bf 1}+{\bf 1}$, so by additivity of the 
monad, $T(2)=T({\bf 1})\times T({\bf 1})=R\times R$); similarly, $E$ is $\mu 
_{{\bf 1}}$, and the diagram commutes by naturality of $\mu$ 
(recalling that $+$ is $T(2\to {\bf 1})$).

\section{Differential calculus for distributions}\label{ad19x}

We attempt in this Section to show how some differential 
calculus of distributions  may be developed independently of the
standard differential calculus of functions. 

For simplicity, we only consider functions and distributions on $R$ 
(one-variable calculus), but some of the considerations  readily generalize 
to distributions and functions on any space $X$ equipped with a 
vector field.

For this, we assume that the monad $T$ on 
$\E$ has the properties described at the end of Section \ref{ad16x}, so in 
particular, $R$ is a commutative ring. To have some differential 
calculus going for such $R$, one needs some further assumption: 
either the possibility to form ``limits'' in $R$-modules, or else, 
availability of the method of synthetic differential geometry (SDG). We 
present the considerations in the latter form, with (some version of) the 
KL\footnote{Kock-Lawvere} axiomatics as the base, since this is well 
suited to make sense in any Cartesian closed category $\E$. In the 
spirit of SDG, we shall talk about the objects of $\E$ as if they 
were sets.

Consider  a commutative ring $R \in \E$. Let $D\subseteq R$ be a 
subset containing 0, and let $V$ be an $R$-module. We say that $V$ 
satisfies the  ``KL''-property (relative to $D$)'' if for any $X\in \E$ 
and any $F:X\times D \to V$, there exists a unique $f':X\to V$ such 
that
$$F(x,d) = F(x,0)+ d\cdot f'(x)\mbox{\quad for all $d\in D$}.$$
If $V$ is a KL module, then so is $X\p V$.

\medskip

Example: 1) models $R$ of synthetic differential geometry, (so  $D$ is the set of $d\in R$ with 
$d^{2}=0$); then the ``classical'' KL axiom says that at least the 
module $R$ itself  satisfies the above condition.
If $X=R$ and $F(x,d)= f(x+d)$ for some function $f:R\to V$, $f'$ is 
the standard derivative of $f$.

2) Any commutative ring, with $D=\{0, d\}$, for one single 
invertible $d\in R$. In this case, 
for given $F$, the $f'$ asserted by the axiom is the function
\begin{equation*}f'(x)= \frac{1}{d}\cdot (F(x,d)-F(x,0));\end{equation*}
if $X=R$ and $F(x,d)= f(x+d)$ for some function $f:R\to V$, $f'$ is 
the standard difference quotient.

\medskip

In either case, we may call $f'$ the {\em derivative} of $f$.

It is easy to see that {\em any} commutative ring $R$ is a model, 
using $\{0, d\}$ as $D$, as in Example 2) (and then also, any $R$-module $V$ 
satisfies then the axiom); this leads to some calculus 
of finite differences. Also, it is true that if $\E$ is the category 
of abstract sets, there are {\em no} non-trivial models of the type in Example 1); 
but, on the other hand, there are other cartesian closed categories $\E$ 
(e.g.\ certain toposes containing the category of smooth manifolds, 
cf.\ e.g.\ \cite{SDG}), and where  
a rather full fledged differential calculus for functions 
emerges from the KL-axiom.

We assume that $R=T(1)$ has the KL property  (for 
some fixed $D\subseteq R$), more generally, that any $R$-module of the form 
$T(X)$ has it.

\begin{prop}[Cancelling universally quantified $d$s]\label{cuqdx} If $V$ is an 
$R$-module which satisfies KL, and $v\in V$ has the property that 
$d\cdot v =0$ for all $d\in D$, then $v=0$.
\end{prop}
{\bf Proof.} Consider the function $f:R\to V$ given by $t\mapsto 
t\cdot v$.
Then for all $x\in R$ and $d\in D$
$$f(x+d)= (x+d)\cdot v = x\cdot v + d\cdot v,$$ 
so that the constant function with value $v$ will serve as $f'$. On 
the other hand, $d\cdot v=d\cdot 0$ by assumption, so that the equation may be continued,
$$= x\cdot  v +d\cdot 0 $$
so that the constant function with value $0\in V$ will likewise serve 
as $f'$. From the uniqueness of $f'$, as requested by the axiom, then 
follows that $v=0$.

\medskip

 We are 
now going to provide a notion of {\em derivative} $P'$ for any $P\in 
T(R)$. Unlike differentiation of distributions in the sense of Schwartz, 
which is defined in terms of differentiation of test functions $\phi$, 
 our construction does not 
mention test functions, and the Schwartz definition\footnote{to be 
precise, Schwartz includes a minus sign in the definition, to 
accommodate the viewpoint that ``{\em distributions are generalized 
functions}'', so that a minus sign in the definition is needed to have 
differentiation of distributions {\em extending} that of functions.}
$\langle P', \phi \rangle := \langle P, \phi '\rangle$ comes in 
our treatment out as a {\em result}, see Proposition \ref{switchy} below.

For $u\in R$, we let $\alpha^{u}$ be ``translation by $u$'', i.e.\ 
the map\footnote{Note that this map is invertible. So we are not 
using the functor property of $T$, except for invertible maps. This 
means that some of the following readily makes sense for distribution 
notions $T$ which are only functorial for a more restricted class of 
maps, like proper maps.} $x\mapsto x+u$. 

For $u=0$,  $\alpha^{u}_{*}(P)-P =0 
\in T(R)$. 
 Assuming that the $R$-module $T(R)$ has the KL property, we therefore have  for 
any $P\in T(R)$ that there exists a unique $P' \in T(R)$ such that for 
all $d\in D$,
$$d\cdot P' =  \alpha^{d}_{*}(P)-P.$$
Since $d\cdot P'$ has total 0 for all $d\in D$, it follows from 
Proposition \ref{cuqdx} that $P'$ 
has total 0.

(If $V$ is a KL module, differentiation of functions $R\to V$ can be 
likewise be described ``functorially'' as a map $R\p V \to R\p V$, 
namely to $\phi \in R\p V$, $\phi'$ is the unique element in $R\p V$ 
satisfying
$d\cdot \phi'= (\alpha^{d})^{*}(\phi ) - \phi$.)

  Differentiation is translation-invariant:
using
$$\alpha ^{t}\circ \alpha ^{s}= \alpha ^{t+s}=\alpha ^{s}\circ \alpha 
^{t},$$
it is easy to deduce that
\begin{equation}\label{transdif}
(\alpha ^{t}_{*}(P))' = (\alpha ^{t})_{*}(P').
\end{equation}
\begin{prop}Differentiation of 
distributions on $R$ is a $T$-linear process.
\end{prop}
{\bf Proof.} Let temporarily $\Delta : T(R)\to T(R)$ denote the differentiation 
process. Consider a fixed $d\in D$. Then for any $P\in T(R)$,
$d\cdot \Delta (P) =d\cdot P'$ is $\alpha^{d}_{*}P-P$; it is a difference of  two 
$T$-linear maps, namely the identity map on $T(R)$ and 
$\alpha^{d}_{*} =T(\alpha ^{d})$, and as such is $T$-linear. Thus for 
each $d\in D$, the map $d\cdot \Delta : T(R)\to T(R)$ is $T$-linear. 
Now to prove $T$-linearity of $\Delta$ means, by monad theory, to prove equality of two 
maps $T^{2}(R)\to T(R)$; and since $d\cdot \Delta$ is $T$-linear, as 
we proved, it follows that the two desired maps $T^{2}(R)\to T(R)$ 
become equal when post-composed with the map ``multiplication by $d$'': 
$T(R)\to T(R)$. Since $d\in D$ was arbitrary, it follows from the 
principle of cancelling universally quantified $d$s (Proposition 
\ref{cuqdx}) 
 that the two desired maps are equal, 
proving $T$-linearity.

\medskip

\begin{prop}\label{totexp}Let $P\in T(R)$. Then $$  E(P')=\tot (P).$$
\end{prop}
{\bf Proof.} The Proposition say that  two maps $T(R)\to R$ agree, 
namely  $E\circ \Delta$ and $\tot $, where $\Delta$, as above, is the 
differentiation process $P\mapsto P'$. Both these maps are 
$T$-linear, so it suffices to prove that the equation holds for the 
case  $P=\eta_{R} (x)$, for each fixed $x\in R$. Let us write 
$\delta_{x}$ for the distribution $\eta_{X}(x)$, so that we do not 
confuse it with a function. So we should prove
$$ E(\delta_{x}')=\tot (\delta _{x}).$$
 By the principle of cancelling 
universally quantified $d$s (Proposition \ref{cuqdx}), it suffices
to prove that for all $d\in D$ that
\begin{equation}\label{cccx}  d\cdot E(\delta_{x}')= d\cdot \tot 
(\delta _{x}).\end{equation}
  The left hand side of (\ref{cccx}) is 
\begin{align*}E(d\cdot \delta _{x}') &= E( 
\alpha^{d}_{*}\delta_{x}-\delta_{x})\\
&= 
E( \delta_{x+d}-\delta_{x})\\
&=E(\delta _{x+d})- E(\delta _{x}) =  (x+d)-x =d,\end{align*}
by the second equation in
(\ref{expdelx}). The right hand side of (\ref{cccx}) is $d$, by 
the first equation in (\ref{expdelx}).
This proves the Proposition.

\medskip

The differentiation process for functions, as a map $R\p V \to R\p 
V$, is likewise $T$-linear, but this important information cannot be 
used in the same way as we used $T$-linearity of the differentiation
$T(R)\to T(R)$, since, unlike $T(R)$, $R\p V$ (not even $R\p R$)  is not known
 to be freely generated by elementary quantities 
like the $\delta _{x}$s. 

\medskip

Here is an important relationship between differentiation of 
distributions on $R$, and of functions $\phi : R\to T(X)$; 
such functions can be differentiated, since $T(X)$ is assumed to be KL  as an 
$R$-module. (In the Schwartz theory, this relationship, with $X={\bf 
1}$, serves as {\em 
definition} of derivative of distributions, except for a sign change, 
see an earlier footnote.)
\begin{prop}\label{switchy}For $P\in T(R)$ and $\phi \in R\p T(X)$, one has
$$\langle P' ,\phi \rangle = \langle P, \phi ' \rangle.$$
\end{prop}
{\bf Proof.} We are comparing two maps $T(R)\times (R\p T(X)) \to 
T(X)$, both of which are $T$-linear in the first variable. Therefore, 
it suffices to prove the equality for the case of $P=\delta _{t} 
(=\eta_{R}(t))$; in fact, by $R$-bilinearity of the pairing and 
Proposition \ref{cuqdx}, it suffices to prove that for any $t\in R$ and 
$d\in D$, we have
$$\langle d\cdot (\delta_{t})', \phi \rangle = \langle \delta _{t}, d\cdot \phi'\rangle.$$ 
The left hand side is $\langle  
\alpha^{d}_{*}(\delta_{t})-\delta_{t}, \phi \rangle$, and using bi-additivity of 
the pairing and (\ref{extra-x}), this gives
$( (\alpha^{d})^{*})(\phi )(t)-\phi (t) =  \phi (t+d)-\phi (t) $, which 
is $d\cdot \phi '(t)$.

\medskip

It is easy to see that if $F:V\to W$ is an $R$-linear map between KL 
modules, we have
\begin{equation}
\label{lindiffx} F\circ \phi' =(F\circ \phi)'
\end{equation}
for any $\phi :R\to V$.

\begin{prop}Let $P\in T(R)$ and $Q\in T(R)$. Then
$$(P*Q)'= P'*Q= P*Q'.$$
\end{prop}{\bf Proof.} Let us  prove 
that $(P*Q)'= P'*Q$ (then $(P*Q)=P*Q' $ follows by commutativity of 
convolution). Both sides depend in a $T$-bilinear way on $P$ 
and $Q$, so it suffices to see the validity for the case where 
$P=\delta_{a}$ and $Q=\delta _{b}$. To prove 
$(\delta_{a}*\delta_{b})'=\delta_{a}'*\delta_{b}$, it suffices to 
prove that for all $d\in D$, 
$$d\cdot (\delta_{a}*\delta_{b})'=d\cdot \delta_{a}'*\delta_{b},$$
and both sides come out as
$\delta _{a+b+d}-\delta_{a+b}$, using that $*$ is $R$-bilinear.

\subsection*{Primitives of distributions on $R$}\label{PEQX}

We noted already in Section \ref{ad11x} that $P$ and $f_{*}(P)$ have 
same total, for any $P\in T(X)$ and $f:X\to Y$. In particular, for 
$P\in T(R)$ and $d\in D$, $d\cdot P' = P-\alpha^{d}_{*}(P)$ has total 
$0$, so cancelling the universally quantified $d$ we get  that $P'$ has 
total $0$.

A {\em primitive} of a distribution $Q\in T(R)$ is a $P\in 
T(R)$ with $P'=Q$. Since any $P'$ has total 0,
a necessary condition that a distribution $Q\in T(R)$ has a primitive is that $\tot 
(Q)=0$. Recall that primitives, in ordinary 1-variable calculus,  are 
also called 
``indefinite {\em integrals}'', whence the following use of the word 
``integration'':

\medskip
\noindent{\bf Integration Axiom.} {\em Every $Q\in T(R)$ with $\tot (Q)=0$ has a unique 
primitive.}

\medskip
(For contrast: for functions $\phi  
:R\to R$, the standard integration axiom is that 
primitives {\em always} exist, but are {\em not} unique, only up to an additive 
constant.) 

By $R$-linearity of the differentiation process $T(R)\to T(R)$, the 
uniqueness assertion in the Axiom is equivalent to the assertion: 
{\em if $P'=0$, then $P=0$}. (Note that $P'=0$ implies that $P$ is 
invariant under translations $\alpha^{d}_{*}(P)=P$ for all $d\in D$.) 
The reasonableness of this latter 
assertion is a two-stage argument: 1) if $P'=0$, $P$ is invariant 
under {\em arbitary translations}, $\alpha ^{u}_{*}(P)=P$. 2) if $P$ 
is invariant under all translations, and has compact support, it must 
be 0. (Implicitly here is: $R$ itself is not compact.)

In standard distribution theory, the Dirac distribution $\delta_{a}$ 
(i.e.\ $\eta_{R}(a)$) 
(where $a\in R$) has a primitive, namely the Heaviside ``function''; but 
this ``function'' has not 
compact support - its support is a half line $\subseteq R$.

On the other hand, the integration axiom provides a (unique) primitive 
for a distribution of the form $\delta_{b}-\delta_{a}$, with $a$ and 
$b$ in $R$. This primitive is denoted $[a,b]$, the ``interval'' from 
$a$ to $b$; thus, the defining equation for this interval is
$$[a,b]'=\delta_{b}-\delta_{a}.$$
Note that 
the phrase ``interval from \ldots to \ldots " does not imply that we 
are considering an ordering $\leq$ on $R$ (although ultimately, one 
wants to do so).
\begin{prop} The total of $[a,b]$ is $b-a$. 
\end{prop} 
{\bf Proof.} We have
\begin{align*} \tot ([a,b]) &=E([a,b]')= E(\delta _{b}-\delta _{a})\\
\intertext{by Proposition \ref{totexp} and the fact that $[a,b]$ is 
a primitive of $\delta _{b}-\delta _{a}$}
&=E(\delta_{b})-E(\delta _{a}) = b-a,
\end{align*}
by (\ref{expdelx}).

\medskip

 It is of some interest to study the sequence of distributions
$$[-a,a], \quad [-a,a]*[-a,a], \quad  [-a,a]*[-a,a]*[-a,a], \quad \ldots ;$$
they have totals $2a, (2a)^{2},(2a)^{3},\ldots$; in particular, if 
$2a=1$, this is a sequence of probability distributions, approaching 
a Gauss normal distribution (the latter, however, has presently no 
place in our context, since it does not have compact support).

\medskip

The following depends on the Leibniz rule for differentiating a 
product of two functions; so this is {\em not} valid under he 
general assumptions of this Section, but needs the further assumption 
of Example 2, namely that $D$ consists of $d\in R$ with $d^{2}=0$, 
as in synthetic differential geometry. 
We shall then use ``test function'' technique to prove a ``Leibniz 
rule''\footnote{with the Schwartz convention of sign, one gets a plus 
rather than a minus between the terms.} for the action $\vdash$ 
\begin{prop}\label{leibnizx}For any $P\in 
T(R)$ and $\phi \in R\p R$,
$$(P\vdash \phi )' = P'\vdash \phi - P\vdash \phi '.$$
\end{prop}
{\bf Proof.} Since there are enough $B$-valued test functions 
(Proposition \ref{testx}), it suffices to prove for arbitrary test 
function $\eta :R \to B$ (with $B$ a free $T$-algebra - hence, by 
assumption, a KL-module) we have
$$(\langle P\vdash \phi )' , \eta \rangle =\langle P'\vdash \phi 
-P\vdash \phi ', \eta 
\rangle.$$
We calculate (using  that the pairing is bi-additive):
\begin{align*}\langle (P\vdash \phi )' , \eta \rangle &= \langle 
P\vdash \phi , \eta '\rangle \mbox{\quad (by Proposition \ref{switchy})}\\
&= \langle P,\phi \cdot \eta' \rangle \mbox{\quad (by Proposition 
\ref{switchx}}\\
&= \langle P, (\phi\cdot \eta)'-\phi'\cdot \eta \rangle\\
\intertext{using that Leibniz rule applies to any bilinear pairing, 
like the multiplication map $\cdot$,}
&=\langle P, (\phi\cdot \eta)'\rangle -\langle P,\phi'\cdot \eta 
\rangle \\
&=\langle P', \phi \cdot \eta \rangle - \langle P, \phi' \cdot \eta 
\rangle\\
\intertext{using Proposition \ref{switchy} on the first summand}
&=\langle P'\vdash \phi ,\eta \rangle - \langle P\vdash \phi' ,\eta 
\rangle \\
\intertext{using Proposition 
\ref{switchx}
 on each summand}
&=\langle P'\vdash \phi - P\vdash \phi' ,\eta \rangle
\end{align*}
  In other words,  the proof looks 
formally like the one from books on distribution theory, but does not 
depend on ``sufficiently many test functions with values in $R$''.

\noindent \verb+kock@imf.au.dk +

\noindent August 2011

\end{document}